\documentclass[12pt]{article}

\usepackage{amsmath,amsthm}
\usepackage{graphicx}
\usepackage{hyperref}
\usepackage{exscale,relsize}
\usepackage{color}
\usepackage{url}
\usepackage{amsfonts}

\usepackage{multicol}

\def\hksqrt{\mathpalette\DHLhksqrt}
\def\DHLhksqrt#1#2{\setbox0=\hbox{$#1\sqrt{#2\,}$}\dimen0=\ht0
\advance\dimen0-0.2\ht0
\setbox2=\hbox{\vrule height\ht0 depth -\dimen0}%
{\box0\lower0.4pt\box2}}
\theoremstyle{definition}

\def\bee{\begin{equation}}
\def\eee{\end{equation}}
\def\OO{\mathcal{O}}
\def \Li{{\rm Li}}
\def \BB{{\cal B}}
\def \TT{{\cal T}}
\def \II{{\cal I}}
\def \h{\hspace{0.15truecm}}
\def\hksqrt{\mathpalette\DHLhksqrt}
\def\DHLhksqrt#1#2{\setbox0=\hbox{$#1\sqrt{#2\,}$}\dimen0=\ht0
\advance\dimen0-0.2\ht0
\setbox2=\hbox{\vrule height\ht0 depth -\dimen0}%
{\box0\lower0.4pt\box2}}

\allowdisplaybreaks

\makeatletter
\@addtoreset{footnote}{page}
\makeatother

\begin{document}

\thispagestyle{empty}
\bigskip
\centerline{    }
\centerline{\Large\bf  Some heuristics  on the gaps between consecutive primes}
\bigskip\bigskip
\centerline{\large\sl Marek Wolf}
\begin{center}
Cardinal  Stefan  Wyszynski  University, Faculty  of Mathematics and Natural Sciences. \\  
ul. W{\'o}ycickiego 1/3,   PL-01-938 Warsaw,   Poland, e-mail:  m.wolf@uksw.edu.pl
\end{center}

\begin{abstract}
We give a review of the rigorous results on the  gaps between consecutive primes. Next we present heuristic arguments
leading  to the formula for the number of pairs of consecutive primes   $p_n, p_{n+1}<x$ separated by gap $d=p_{n+1}-p_n$
expressed directly by $\pi(x)$, i.e.  the number of all primes  $<x$. We use  this formula  to discuss the problem of champions,
to find the  maximal gap between two consecutive primes smaller than $x$ represented by $\pi(x)$,
 generalized Brun's constants and next  the new formula for first appearance of primes separated by gap $d$.  We derive  from our guesses the leading term
$\log \log(x)$ in the prime harmonic sum.   Finally  we  discuss the Andrica Conjectures. We  illustrate these topics by
extensive computer data   collected up to $2^{48}=2.81\ldots 10^{14}$.
\end{abstract}

\bigskip\bigskip\bigskip



\bibliographystyle{abbrv}

\section{The overview of the problem and rigorous results}
\label{Introduction.}

To investigate the set of prime numbers $\{2, 3, 5, \ldots, p_n, \ldots\}$ one can follow  many approaches. 
First of all we can ask what is  the number of primes up to a given threshold $x$. This function is usually denoted by $\pi(x)$.
It is one of the greatest surprises  mathematics that such an erratic function
as $\pi(x)$ can be approximated by a simple expression. Namely  Carl Friedrich Gauss as a teenager (different sources put his age
between  fifteen years and seventeen years)   made at the end of the  eighteen century  conjecture that
$\pi(x)$ is roughly  given by the logarithmic integral $\Li(x)$:
\bee
\pi(x)\sim {\Li}(x):=\int_2^x  \frac{du}{{\log(u)}} \approx \frac{x}{\log(x)}.
\label{PNT}
\eee
The   symbol $f(x)\sim  g(x)$ means  that $\lim_{x\rightarrow \infty} f(x)/g(x)=1$.  Presently many proofs of \eqref{PNT}  are
known and this formula is called the  Prime Number Theorem (PNT).
On the other hand  one can look at the differences between arbitrary primes   $d=p'-p$ or at the  distances between consecutive primes
$d_n=p_{n+1}-p_n$.   In 1922   G. H. Hardy and J.E. Littlewood in the famous paper  \cite{Hardy_and_Littlewood}
proposed 15 conjectures. The  conjecture B of their paper states:

{\it
There are infinitely many prime  pairs $(p, p^\prime)$, where $p^\prime = p+d$,
for every even $d$. If $\pi_d(x)$ denotes the number of prime pairs differing
by $d$ and  less than $x$, then}
\bee
\pi_d(x)\sim  C_2  \prod_{p\mid d} \frac{p-1}{p-2}~ \frac{x}{\log^2(x)}.
\label{H-L}
\eee
Here the product is over odd primes  $p\geq 3$ dividing $d$. The twin primes constant $C_2\equiv 2c_2$ is defined by the infinite product:
\bee
C_2 \equiv  2c_2 \equiv 2\prod_{p > 2} \biggl( 1 - {\frac{1}{ (p - 1)^2}}\biggr) =1.32032363169\ldots
\label{stalac2}
\eee

\begin{figure}[ht]
\begin{center}
\includegraphics[width=0.8\textwidth, angle=0]{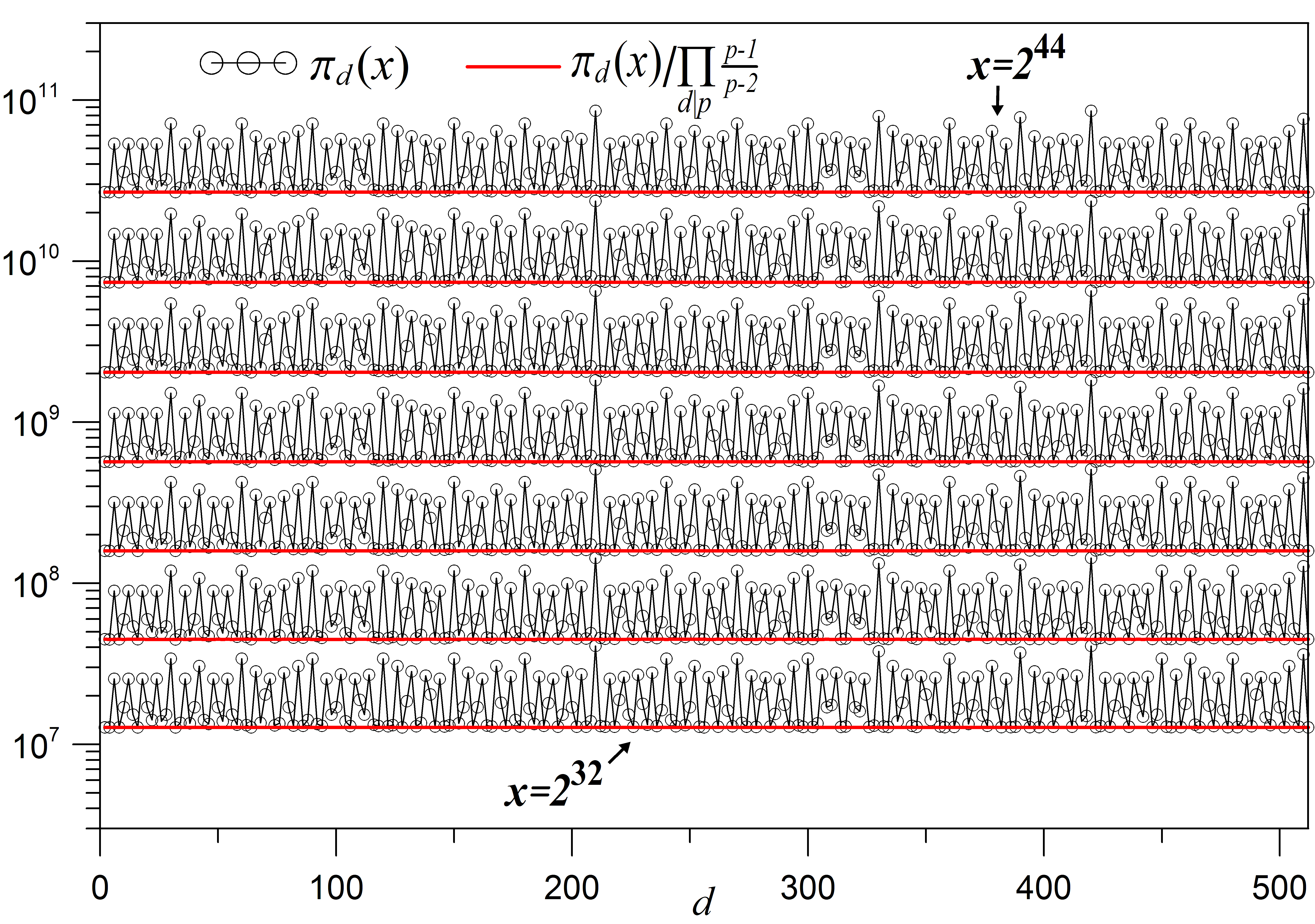} \\
\caption{ The plot of $\pi_d(x)$ (eq. (\ref{H-L})) obtained from the computer search
for $d=2,4, \ldots, 512$ and for $x=2^{32}, 2^{34}, \ldots, 2^{44}$. In red are
the ratios $\pi_d(x)/{\mathfrak S}(d)$ plotted showing explicitly that a characteristic
oscillating pattern with peaks at $6k$, $30k, 210k$ is caused by the product
${\mathfrak S}(d)$.}
\label{fig-H-L}
\end{center}
\end{figure}

Computer results of the search for pairs of primes separated by a
distance $d\leq 512$ and smaller than $x$ for $x=2^{32}, 2^{34}, \ldots,
2^{44}\approx 1.76 \times 10^{13}$  are shown in Fig.\ref{fig-H-L}
and they provide a firm support in favor of (\ref{H-L}).
The  characteristic oscillating pattern of points is caused by the product
\bee
{\mathfrak S}(d)= \prod_{p \mid d, p > 2} \frac{p - 1}{p - 2}
\label{product}
\eee
appearing in (\ref{H-L}).  This product  ${\mathfrak S}(d)$ has local maxima for $d$ equal to  the products of consecutive primes
(factorials over primes $p_n\sharp :=2\cdot 3 \cdot 5 \cdot \ldots \cdot p_n$  are  called ``primorials''):
${\mathfrak S}(6)=2, ~{\mathfrak S}(30)=8/3=2.666\ldots ~{\mathfrak S}(210)=16/5=3.2, \dots$  (local minima are 1 and they correspond to $d=2^m$).
Clearly visible  in Fig. \ref{fig-H-L} are oscillations of the period
$6=2\times 3$ with overimposed  higher harmonics $30=2\times 3\times 5$  and  $210=2\times 3\times 5\times 7$, i.e. when ${\mathfrak S}(d)$
has local maxima.  The red lines present  $\pi_d(x)/{\mathfrak S}(d)$ and they are  perfect straight  lines $C_2 x/\log^2(x)$.
There is large evidence  both analytical and experimental in favor of (\ref{H-L}).
Besides the original circle method used by Hardy and Littlewood
\cite{Hardy_and_Littlewood}  there appeared
papers  \cite{Polya1959} and  \cite{Rubinstein} where other heuristic arguments
were presented.

The above notation $\pi_d(x)$ denotes prime pairs not necessarily successive.
Not much is known about gaps  between {\it consecutive} primes, which  seems to be more interesting and
difficult  than the case of pairs of arbitrary (not consecutive) primes treated by the
Hardy--Littlewood  conjecture B.  The Polignac's conjecture \cite[p. 187]{Ribenboim} asserts that for  every even
natural number $d$  there exist infinitely many consecutive primes $p_n, p_{n+1}$, such that $p_{n+1}-p_n=d$.
Let $\tau_d(x)$ denote the  number of pairs  of {\it consecutive} primes smaller than a given bound $x$ and separated by $d$:
\bee
\tau_d(x)= \{{\rm number~ of~ pairs~~} p_n,  p_{n+1} < x,~~{\rm  with}~
d=p_{n+1}-p_n\}.
\label{definition}
\eee
For odd $d$ we  supplement this definition by putting $\tau_{2k+1}(x)=0$ (although $\tau_1(x)=1$ for all $x\geq 3$).
The pairs of primes separated by $d=2$ (twins)  and $d=4$  (cousins)   are special as they
always have to be consecutive primes (with the exception of the pair
(3,7) containing 5 in the middle):  $\pi_2(x)\equiv \tau_2(x), ~\pi_4(x)\equiv \tau_4(x)+1$  and the  Hardy--Littlewood conjecture B
gives that there is approximately  the same number of twins and cousins:   $\pi_2(x)\approx  \pi_4(x) \sim C_2 \frac{x}{\log^2(x)}$.
In this paper we will present simple
heuristic reasoning leading to the formula for   $\tau_d(x)$ {\it expressed
directly by} $\pi(x)$ --- the total number of primes up to $x$.

A few  main questions related to  the problem of gaps $d_n = p_{n+1} - p_n$
between consecutive primes can be distinguished. From  $\pi(x)\sim x/\log(x)$ it follows that the mean gap between
consecutive primes  is  of the order $\log(x)$.   There are gaps of arbitrary length
between primes: namely the $n$ numbers $(n+1)!+2, (n+1)!+3, (n+1)!+4, \dots,
(n+1)!+n+1$ are all composite.  In fact a gap of size $d$  appears much earlier than at $(d+1)!$, see  Section  \ref{first-d}.  The Bertrand's  postulate,  that there is
always a prime  between  $n$  and  $2n$,   was  proved by Chebyshev  in 1852.  The Bertrand's postulat in another formulation
says that $d_n<p_n$ for every $n\geq 1$.
The growth rate of the form $d_n = \OO(p_n^\theta)$ with different $\theta<1$ \footnote{The big--$\OO$ symbol
$ f(x) = \OO (g(x))$  as  $x \to \infty $ for two functions $f(x)$
and  $g(x)$   means that  there exists   positive constant $M$ such that for all
sufficiently large values of $x$ the inequality   $|f(x)|\leq \;M|g(x)|$  holds.}  was  proved in the past.
A  few results with $\theta$ closest to 1/2 are the  results of:  C. Mozzochi
\cite{Mozzochi1986}  $\theta = \frac{1051}{1920}$, S. Lou and Q. Yao  obtained
$\theta=6/11$ \cite{Lou1992}, R.C.  Baker and  G. Harman have improved it
to $\theta = 0.535$ \cite{Baker1996} and recently R.C. Baker
G. Harman and J. Pintz  \cite{Baker-et-al} have improved it by 0.01  to
$\theta=21/40=0.525$ which currently  remains the best unconditional result.
The Riemann Hypothesis implies $d_n=\OO(\sqrt{p_n}\log(p_n))$ and $\theta = \frac{1}{2} + \epsilon$ for any $\epsilon>0$.
For a review of results on $\theta$ see \cite{Pintz-Landau}.  Another kind  of results concerns so called small gaps between
primes, i.e.  gaps smaller than the mean gap $\log(x)$.  One can  compare $d_n=p_{n+1}-p_n$  with $\log(p_n)$ and look what
is the limes inferior of the sequence $d_n/\log(p_n)$.  There  were many estimations of this limit culminating with the famous
theorem GPY1   \cite{GPYI}:
\[
\liminf_{n\to \infty} \frac{p_{n+1}-p_n}{\log p_n} =0.
\]
The above formulas for the speed of the growth of $d_n$ are connected with the Legendre's conjecture that there is a prime number
between $n^2$ and $(n + 1)^2$ for every positive integer $n$.  Indeed,  the distance between  $n^2$ and $(n + 1)^2$ is $2n+1$  and
if the  gap  between  primes at $x^2=n^2$ is of the order $\OO(x^{2\theta}) =Cx^{2\theta}$ then  to   have at least one prime between
$x^2$  and  $(x+1)^2$  one has to require $Cx^{2\theta} < 2x $  and  hence  the proof of $\theta<\frac{1}{2}$  is  needed.
In  fact  there is  usually a lot  of  primes  between $n^2$ and $(n + 1)^2$, see \href{https://oeis.org/A014085}{OEISA014085}.
From Gauss's  formula  \eqref{PNT} we have $\pi((n+1)^2)-\pi(n^2) \sim  n/\log(n)$.  Let us remark that the fact $\theta < 2/3$
suffices to show that between $n^3$ and $(n+1)^3$  there is always  at least one prime:   the gap  between two consecutive
primes around $n^3$ is $\OO(n^{3\theta})$  and  for $\theta<2/3$  it is smaller  than distance  from  $n^3$ to  $(n+1)^3$
which is $\OO(n^2)$.

In the middle of 2013 the major step towards the proof of the conjecture B and Polignac's conjecture  was made: Yitang Zhang
published in  Annals of Mathematics the paper \cite{Zhang-2014}  in which he proved unconditionally that
there exists gap $d<7 \times 10^7$ which is a difference of infinitely many pairs of two primes and that
$\liminf_{n\to \infty}\, (p_{n+1} - p_n) < 7 \times 10^7$. It means,  that there is, at least one,  such $d<7 \times 10^7$ that
there exist infinity of primes pairs separated by $d$.   This achievement brought to Zhang great fame and popularity: there were
many articles in daily and weekly press, see e.g. \cite{Zhang_new-yorker}.   Very soon his bound $7 \times 10^7$  was lowered many
orders  by  mathematicians and the separate projects Polymath 8a \cite{Polymath-RIMS} and Polymath 8b  \cite{Polymath}.
J. Maynard \cite{Maynard-2015} proved unconditionally that   $\liminf_{n\to \infty}\, (p_{n+1}-p_n)\le 600$ and current record
is 246  obtained by Polymath 8b.  Assuming the Elliott-Halberstam conjecture Maynard showed that $\liminf_{n\to \infty}\, (p_{n+1} - p_n) < 12$
and Polymath  lowered it to 6.

\begin{figure}
\begin{center}
\vspace{-0.7cm}\hspace{-1.7cm}
\includegraphics[height=0.7\textheight, angle=90]{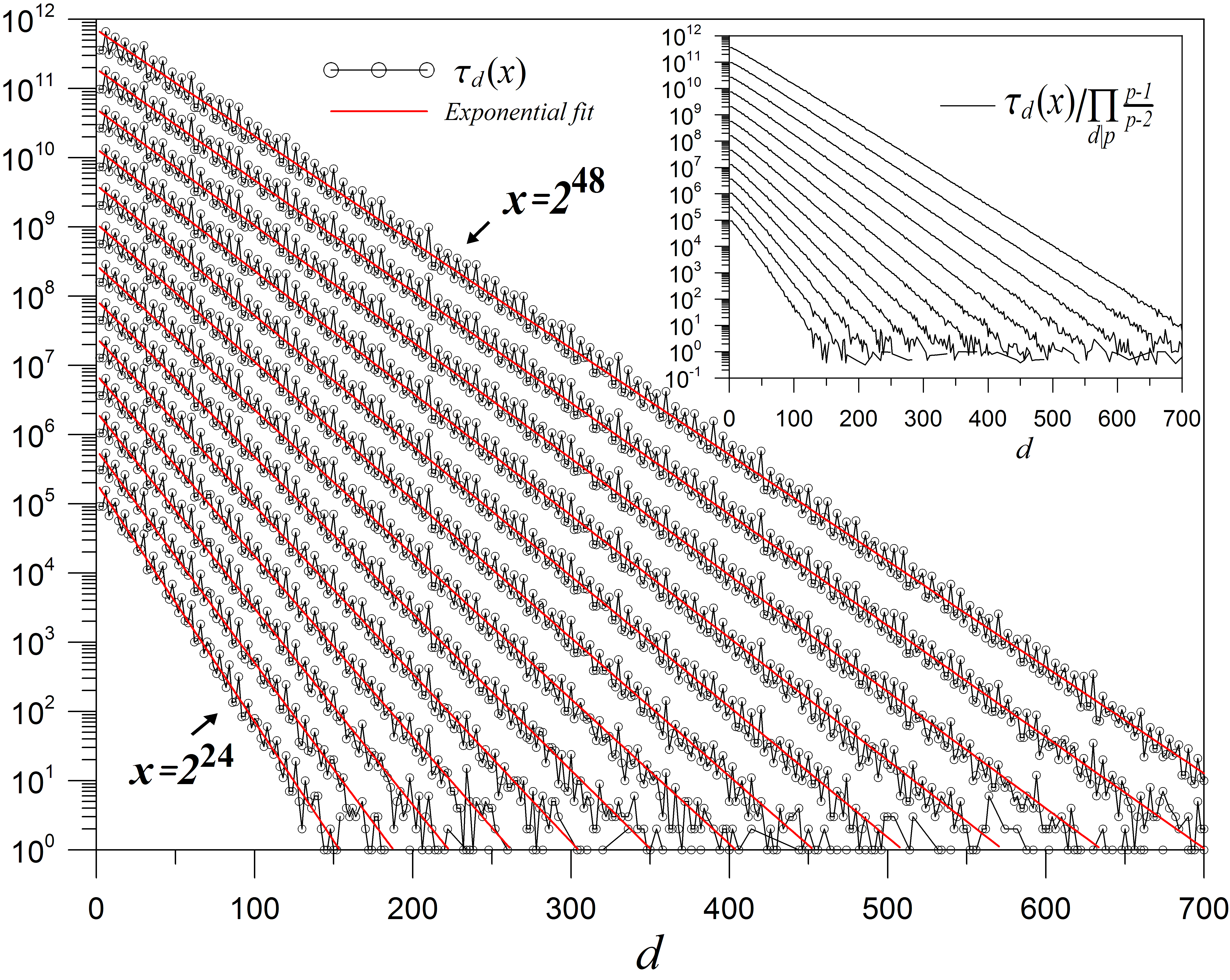} \\
\caption{Plots of $\tau_d(x)$ for $x=2^{24}, 2^{26}, \ldots, 2^{46}, 2^{48}$.
In the inset plots of $\tau_d(x)/{\mathfrak S}(d)$ are shown for the same set of $x$.
In red are exponential fits $a(x)e^{-db(x)}$ plotted. Spikes at  $d=6$, $d=30, d=210$ are responsible for
the champions phenomenon.}
\label{fig-hist}
\end{center}
\end{figure}

On the opposite side  the second question about $d_n$ concerns the existence of very large gaps, i.e.  gaps of the order $\log^2(x)$.
Let  $G(x)$  denotes the largest gap between consecutive primes below a
given bound $x$:
\bee
G(x)={\max_{p_n<x}}~(p_n-p_{n-1}).
\label{def-G}
\eee
For this function lower bounds $f(x)$ are searched for: $G(x)>f(x)$. The Prime Number Theorem
$\pi(x)\sim x/\log(x)$ trivially gives $G(x)>\log(x)$. Better inequality
\bee
G(x) \geq \frac{(c e^\gamma + o(1))\log(x) \log\log(x) \log\log\log\log(x)}{(\log \log \log (x))^{2}},
\eee
where $\gamma=0.577216\ldots$ is the Euler-Mascheroni constant, was proved by
H. Maier and C. Pomerance in \cite{Pomerance1990} with $c=1.31256\ldots$ and improved
by J. Pintz to $c=2$ in \cite{Pintz1997}.  J. Maynard in \cite{Maynard-2016} proved that
\bee
G(x)> \frac{c(1+o(1))\log(x) \log\log(x) \log\log\log\log(x)}{(\log \log \log (x))^{2}}
\label{eq-Maynard}
\eee
for any fixed $c$, i.e. in the results of Maier, Pomerance and Pintz one can take the constant $c$ to be arbitrarily large.
It was independently proved in   \cite{FGKT}  in another form.
Recently in \cite{FGKMT}  the denominator in \eqref{eq-Maynard}  was  lowered (hence  $G(x)$ attains  grater values):
\bee
G(x)> \frac{c(1+o(1))\log(x) \log\log(x) \log\log\log\log(x)}{\log \log \log (x)}.
\eee

In 1946 there appeared a paper \cite{Cherwell1946}, where the problem
of different patterns of pairs, triplets etc. of primes was treated by the
probabilistic methods. In particular the formula for a number of primes$<x$
and separated by a gap $d$ was deduced on p. 57 from probabilistic arguments.
Similar equation  appeared many years later  in the paper by S. Guiasu in  \cite{Guiasu-1995}. We will discuss these   formulas  in next Section.

In 1974 there appeared a paper by Brent \cite{Brent1974}, where
statistical properties of the distribution of gaps between
consecutive primes were studied both theoretically and numerically. Brent
had  applied the inclusion--exclusion principle and obtained from
(\ref{H-L}) a formula for the number of consecutive prime pairs less
than $x$ separated by $d$. But
his result (formula (4) in \cite{Brent1974}) does  not have a closed form
and he had to produce   on a  computer   the  table of constants appearing in his formula (4). The
attempt to estimate these sums and to write a closed formula for them
was undertaken  in \cite{champions}. Further development of ideas contained in this paper was published by  D. A. Goldston
and A. H. Ledoan \cite{Goldston-Leodan-2012} in 2012.
Below we will present  heuristic considerations leading to closed formulas for  some quantities characterizing gaps
between consecutive primes.


\section{Heuristic   formula  for $\tau_d(x)$}    
\label{main-conjecture}

\indent   During over a seven months long run of  the  program  on the  64--bits 2.7 GHz  computer
we have collected   values of $\tau_d(x)$ up to $x=2^{48}\approx 2.8147\times 10^{14}$.
The data representing the function $\tau_d(x)$ were stored at values of $x$ forming the geometrical
progression with the ratio 2, i.e. at $x=2^{15}, 2^{16}, \ldots, 2^{47},
2^{48}$;  the largest encountered  gap was $d=906$ . Such a choice of the intermediate thresholds as powers of 2
was determined by the employed computer program in which  the primes were coded
as bits. The resulting curves are plotted in Fig.\ref{fig-hist}.
The data is available for downloading from  \url{http://pracownicy.uksw.edu.pl/mwolf/gaps.zip}.

In the plots of $\tau_d(x)$ in  Fig.\ref{fig-hist} a lot of regularities can
be observed.
The pattern of points in Fig.\ref{fig-hist} {\it does not
depend on $x$}: for each $x$ the arrangements of circles is the same,
only the intercept increases and the slope decreases. Like in the case of $\pi_d(x)$
the oscillations are described by the product ${\mathfrak S}(d)$, see the inset in Fig. \ref{fig-hist}.
The fact that the  points in Fig.\ref{fig-hist} lie around the straight lines on the
semi-logarithmic scale suggest for $\tau_d(x)$ the following

{ \bf Ansatz 1 :}
\bee
\tau_d(x)={\mathfrak S}(d) B(x) F^{d}(x),
\label{ansatz}
\eee
where $F(x)<1$ (because $\tau_d(x)$ decreases with $d$).

The essential point of the presented below considerations consists in a possibility
of determining the two unknown functions  $B(x)$  and  $F(x)$  by
{\it assuming only the above exponential decrease of $\tau_d(x)$ with $d$ and
employing two identities fulfilled by} $\tau_d(x)$  {\it just by definition}.
First of all, the number of all gaps is  equal to  the number of all    
primes smaller than $x$  minus 2:
\bee
\sum_{ (p_{n}-p_{n-1} ),~p_{n}<x} 1 \equiv \sum_{d=2}^{G(x)} \tau_d(x)~=~~\pi(x) - 2,
\label{identity_1}
\eee
where $G(x)$ is the largest gap below $x$  which was defined in (\ref{def-G}).
The second self--consistency condition comes from an observation  that the
sum of differences between consecutive primes $p_n \leq x$ is equal to the
largest prime $\leq x$ (minus 3 coming from the distance to $p_2=3$)
and for large $x$ we can write:
\bee
\sum_{p_{n}<x}  (p_{n}-p_{n-1} )\equiv \sum_{d=2}^{G(x)} d~\tau_d(x) \approx ~x.
\label{identity_2}
\eee
The erratic behavior of the product ${\mathfrak S}(d)$ is an obstacle in
calculation of the above sums (\ref{identity_1}) and  (\ref{identity_2}).
We will replace the product ${\mathfrak S}(d)$ in the sums by its average value.
In \cite{Bombieri} E. Bombieri and H. Davenport  have proved that:
\bee
\sum_{k=1}^n \prod_{p\mid k,p>2}\frac{p-1}{ p-2} =
{n \over \prod_{p > 2}( 1 - \frac{1}{(p - 1)^2})} + \mathcal{O}(\log^2(n));
\eee
i.e. in the limit $n\to \infty$ the number  $1/\prod_{p > 2}( 1 - \frac{1}{(p - 1)^2})$
is the arithmetical average of the product $\prod_{p\mid k} \frac{p-1}{p-2}$.
Thus we will assume that for functions $f(k)$ going to zero like $const ^{-k}$ the
following identity holds:
\bee
\sum_{k=1}^\infty \prod_{p\mid k, p>2}\frac{p-1}{p-2} ~f(k) = \frac{1}{\prod_{p
> 2}( 1 - \frac{1}{(p - 1)^2})} \sum_{k=1}^\infty f(k).
\label{rownosc1}
\eee

We can justify the above formula by invoking  the Abel  partial summation  in the form  \cite[Th. 421]{H-W}:
$$
\sum_{k=1}^n a_k b_k = -\sum_{k=1}^{n-1} S(k)c_k + S(n) b_n,
$$
where  $S(k)=a_1+ \dots + a_k$  and $c_k=b_{k+1}-b_k$.
Putting here $a_k={\mathfrak S}(k)$, $b_k=f(k)$, $ S(k)=k/c_2+\mathcal{O}(\log^2(k))$,
next replacing $\log^2(2) < \log^2(3) < \ldots < \log^2(n-1)$ by larger $\log^2(n)$
and collecting terms we obtain in the part multiplied by $1/c_2$ the sum
$f(1)+f(2)+ \ldots +f(n)$ and in the part multiplied by  $\mathcal{O}(\log^2(n))$
we see that the values $f(2), \ldots f(n-1)$ cancel pairwise leaving only $f(1)=0$
and $f(n)$:
\bee
\sum_{k=1}^n {\mathfrak S} (k) f(k) = \frac{1}{c_2}\sum_{k=1}^{n} f(k) + 
f(n)\Big(\frac{n}{c_2}+\mathcal{O}\big(\log^2(n)\big)\Big).
\label{rownosc2}
\eee
Taking the limit $n \to \infty$ we get (\ref{rownosc1}) as the series $\sum_{k=1}^{\infty} f(k)$
converges if $f(k)$  goes to zero sufficiently fast. Thus in the sums of ${\mathfrak S}(k)f(k)$ the product ${\mathfrak S}(k)$
can be replaced by its mean  value $1/c_2=2/C_2=1.51478\ldots$.

We extend in (\ref{identity_1})  and eq. (\ref{identity_2}) the summations to infinity
and using the Ansatz (\ref{ansatz}) and  eq. (\ref{rownosc1}) we get the geometrical and
differentiated geometrical series. For odd $d$ we have defined $\tau_{2k+1}(x)=0$.
Then, writing $d=2k$ we obtain:
\bee
\sum_{d=2, 4,6, \ldots}^{\infty} \tau_d(x)~=~\frac{B(x)}{c_2}\sum_{k=1}^{\infty} F^{2k}(x)=
{\frac{2}{c_2} \frac{B(x)F^{2}(x)}{1-F^{2}(x)}}
\eee
\bee
\sum_{d=2, 4,6, \ldots}^{\infty} d\tau_d(x)~=~\frac{2B(x)}{c_2}\sum_{k=1}^{\infty} k F^{2k}(x)=\frac{1}{c_2} \frac{B(x)F^{2}(x)}{(1-F^{2}(x))^2}.
\eee
By extending  summations in (\ref{identity_1}) and  (\ref{identity_2}) to infinity $G(x)\to \infty$
we made an  error of the order $\mathcal{O}(F(x)^{G(x)+2})$ in the
first case and an error  $\mathcal{O}(G(x)F(x)^{G(x)+2})$ in the second equation, both
going to zero for $x \to \infty$, because for $x \to \infty$ we have $G(x)\to \infty$.
Thus we obtain two equations:
\bee
\frac{1}{c_2}\frac{B(x)F^{2}(x)}{1-F^{2}(x)}  =
\pi(x),\quad\quad \quad \frac{1}{c_2}\frac{2B(x)F^{2}(x)}{(1-F^{2}(x))^2} = x
\label{rownania}
\eee
of which solutions are
\bee
B(x)=\frac{2c_2\pi^2(x)}{x}\frac{1}{(1-\frac{2\pi(x)}{x})},\quad\quad \quad F^{2}(x)= 1-\frac{2\pi(x)}{x}
\eee
and  a posteriori the inequality $F(x)<1$ holds evidently. Finally, we state the main \\

\indent {\bf Conjecture 1}\\
\indent The function $\tau_d(x)$ is expressed directly by $\pi(x)$:
\bee
\tau_d(x) \sim  C_2 \prod_{p \mid d, p > 2} \frac{p - 1}{p - 2}~~ \frac{\pi^2(x)}{x}\Bigg(1-\frac{2\pi(x)}{x}\Bigg)^{{\frac{ d}{2}-1}} ~~~{\rm for} ~ d\geq 6.
\label{main}
\eee
Similar  formula but written in slightly different form and obtained from probabilistic arguments appeared for the first time
apparently in \cite[p.57]{Cherwell1946}. S. Guiasu obtained the formula  \cite[eq.(7)]{Guiasu-1995}  for probability to find
the gap $2k$  (equal in our notation to $d-2$) among
primes   up to  $n$ from the demand that  the entropy associated  with the probability distribution of these gaps is maximal.
His formula does not contain the product $ {\mathfrak S}(d)$.  For twins ($d=2$) and
cousins ($d=4$) the identities $\tau_{2,4}(x) = \pi_{2,4}(x)$ hold. Because $d$ is even the power of $(1-2\pi(x)/x)$ has a finite
number of terms. The formula (\ref{main}) consists of three terms. The first one depends
only on $d$, the second only on $x$, but the third term depends both on  $d$ and $x$.
In the usual probabilistic approach one should obtain $(1-\frac{\pi(x)}{x})^{{d}-1}$,
see e.g. \cite{Heath-Brown-1988}, \cite[p. 3]{Soundararajan-2007}: to  have a pair of
adjacent  primes separated by $d$ there have to be $d-1$ consecutive composite numbers
in between and probability of such an event is $(1-\pi(x)/x)^{d-1}$; then the term in
front of it comes from the normalization condition.

Although (\ref{main}) is postulated for  $d\geq 6$, we get from it for $d=2$  (and $d=4$):

\indent {\bf Conjecture 2}\\
\bee
\tau_2(x)\equiv \pi_2(x)  \sim C_2  \frac{\pi^2(x)}{x}   
\label{twins-m}
\eee
instead of the usual conjectures
\bee
\pi_2(x) \sim C_2 \frac{x}{\log^2(x)}  
\label{twins-a}
\eee
or
\bee
\pi_2(x)\sim C_2 \int_2^x \frac{du}{\log^2(u)} \equiv C_2 {\Li}_2(x).
\label{twins-b}
\eee
The equation (\ref{twins-m}) expresses the intuitively obvious fact that
the number of twins should be proportional to the square of $\pi(x)$.
Of course (\ref{twins-m}) for $\pi(x)\sim x/\log(x)$ goes into (\ref{twins-a}).
We have checked with the available computer data (see http://sweet.ua.pt/tos/primes.html\#t2)  that (\ref{twins-m}) is better
than (\ref{twins-a}) but worse than (\ref{twins-b}).   Because $\Li_2(x)$ in
(\ref{twins-b}) monotonically increases  while there are local
fluctuations in the  density of
primes and twins, the above formula (\ref{twins-m}) incorporates all irregularities
in the distribution of primes into the formula for the number of twins. Since
both $d=2$ and $d=4$ gaps are necessarily consecutive, we propose the identical
expression (\ref{twins-m}) for $\tau_4(x)\equiv \pi_4(x) \approx \pi_2(x)$,
see \cite{Wolf-RW}.

It is possible to obtain another form of the formula for $\tau_d(x)$, more convenient
for later applications.  Namely, let us represent the function $F(x)$ in the form:
$F(x)=e^{-A(x)}$, i.e. now the {   Ansatz 1  has the form

{\bf Ansatz 1{\rm $'$}}
\bee
\tau_d(x) \sim  B(x) {\mathfrak S}(d) e^{-A(x) d},
\label{ansatz2}
\eee
where $A(x)$ is the slope of the lines plotted in red in
Fig. \ref{fig-hist} and as we  can see $A(x)$ goes to zero for $x \to \infty$.
In  the equations (\ref{rownania})  we use in the nominators the approximation
$e^{-2A(x)}\approx 1-2A(x)$ and in the denominators
$1-e^{-2A(x)}\approx 2A(x)$ for small $A(x)$  and we  obtain

{\bf Conjecture $\bf 1'$}
\bee
\tau_d(x) = C_2 \frac{\pi^2(x)}{x-2\pi(x)} \prod_{p \mid d, p > 2 }\frac{p - 1}{p - 2} e^{- {d \pi(x)/x}} + error~ term(x, d)~~{\rm for} ~ d\geq 6.
\label{main2a}
\eee
For large $x$ we can skip $2\pi(x)$ in comparison with $x$
in the denominator and obtain finally the following pleasant formula:\\

{\bf Conjecture $\bf 1''$}
\bee
\tau_d(x) = C_2 \frac{\pi^2(x)}{x} \prod_{p \mid d, p > 2}\frac{p - 1}{p - 2} e^{- {d \pi(x)/x}} + error~ term(x, d)~~{\rm for} ~ d\geq 6.
\label{main2}
\eee

In equation  (\ref{main2})  the term in the exponent has
a simple interpretation: difference $d$ is divided by the mean gap $x/\pi(x)$
between consecutive primes.
Because for small $u$ an  approximation  $\log(1-u)\approx -u$ holds, we can turn for
large $x$ the conjecture (\ref{main}) to the form of conjecture \eqref{main2}:
\bee
\Bigg(1-\frac{2\pi(x)}{x}\Bigg)^{{\frac{ d}{\ 2}}}=e^{\frac{d}{2}\log\big(1-\frac{2\pi(x)}{x}\big)}
\approx e^{-\frac{ d \pi(x)}{x}}.
\eee

Putting in (\ref{main2}) $\pi(x)\sim x/\log(x)$ and comparing with the original
Hardy--Littlewood conjecture we obtain that
the number $\tau_d(x)$ of {\it successive} primes $(p_n, ~p_{n+1})$ smaller
than $x$ and of the difference $d=p_{n+1}-p_n$ is diminished by the factor
$\exp(-d/\log(x))$  in comparison with the number of {\it all} pairs of primes $(p, p^\prime)$ apart  in
the distance $d=p^\prime-p$:
\bee \tau_d(x) \sim  \pi_d(x) e^{-d/\log(x)}  ~~~~~{\rm for~} d\geq 6.
\label{relation}
\eee
Heuristically, this relation encodes in the  series for $e^{-d/\log(x)}=1-d/\log(x) + (d/\log(x))^2/2!-(d/\log(x))^3/3!+\ldots$
the  inclusion-exclusion principle for obtaining $\tau_d(x)$ from
$\pi_d(x)$. The above relation   \eqref{relation} is confirmed by comparing  the Figures
\ref{fig-H-L} and \ref{fig-hist}. R.P. Brent in \cite{Brent1974}
using the inclusion-exclusion principle  has
obtained from the B conjecture of Hardy and Littlewood
the formula for $\tau_d(x)$,  which agrees very well with computer results. However the
formula of Brent (eq.(4) in the paper \cite{Brent1974}) is not of a closed
form: it contains a double sequence of constants $A_{r,k}$,  which can be
calculated only by a direct use of the computer, what is very time consuming, see
discussion of S. Herzog at the web site
http://zigherzog.net/primes/. R. P. Brent in \cite{Brent1974} in Table 2 compares the
number of actual gaps $d=2,\ldots, 80$ in the interval $(10^6, 10^9)$ with the numbers
predicted from his  formula finding perfect agreement. Analogous method to determine
the values of $\tau_d(x)$ was employed in \cite[see eq.(2-8) and the preceding formula]{champions}.
The formula (2-8) from \cite{champions} adapted to our notation has the form:
\bee
\tau_d(x) \sim C_2 {\mathfrak S}(d) \int_2^x \frac{\exp(-d/\log(u))}{\log^2(u)} du.
\eee
Integrating the above integral once by parts gives a term $xe^{-d/\log(x)}/\log^2(x)$
corresponding to (\ref{main2}) with $\pi(x)\sim x/\log(x)$.  The expression (\ref{main2})  for $\tau_d(x)$ was
proved in slightly different form   
under the assumption  of the conjecture B of Hardy--Littlewood by  D. A. Goldston and A. H. Ledoan \cite{Goldston-Leodan-2012}
in 2012.   They proved that for any positive constant $\lambda$ and $d$ even with $d \sim \lambda \log x$ as $x \to \infty$,
we have
\bee
\tau_d(x) \sim  C_2 \mathfrak{S}(d) \frac{x}{(\log x)^2} e^{-\lambda}.
\eee

\begin{figure}[h]
\begin{minipage}[t]{0.5\textwidth}
\centering
\includegraphics[width=\textwidth]{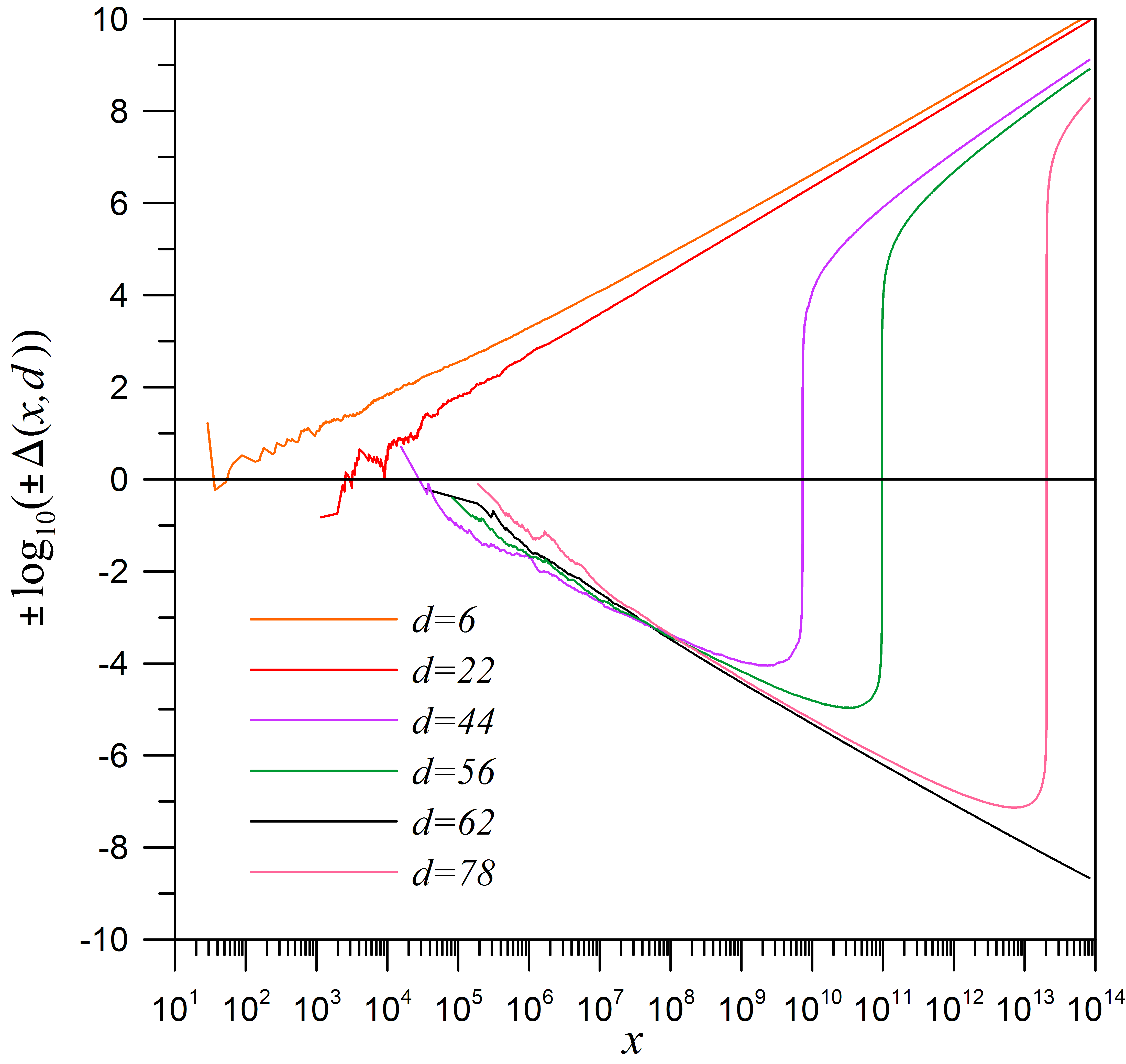}
\caption{Plots of $\Delta(x, d)$ on the double logarithmic scale for $d=6, 22, 44, 56, 62, 78$.
On the $y$ axis we have plotted $\log_{10}(\Delta(x, d))$ if $\Delta(x, d)>0$ and
$-\log_{10}(-\Delta(x, d))$ if $\Delta(x, d)<0$.}
\label{fig-wykres-delty}
\end{minipage}
\hspace{0.03\textwidth}
\begin{minipage}[t]{0.5\textwidth}
\centering
\includegraphics[width=0.95\textwidth]{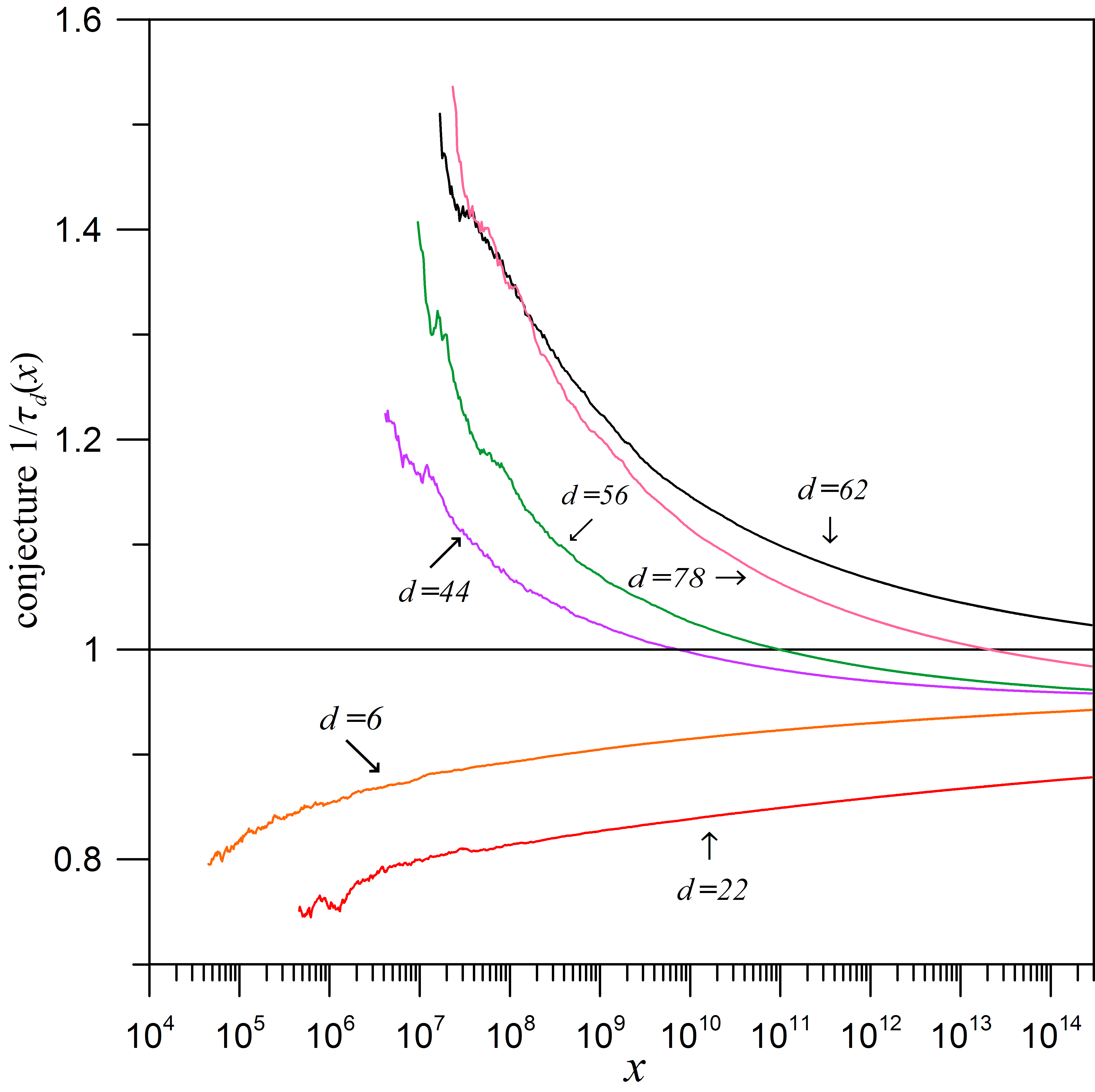}
\caption{Plots of ratios of the values predicted from the
Conjecture 1 to the real values of $\tau_d(x)$ for $d=6, 22, 44, 56, 62, 78$.
The plots begin at such $x$ that $\tau_d(x)>1000$ to avoid large initial fluctuations
of these ratios (see initial parts of curves in the previous Figure).}
\label{fig-compare_tau_errors}
\end{minipage}
\end{figure}

It is not possible to guess an analytical form of error terms in formulas
(\ref{main})  and (\ref{main2}) at present (let us remark that
the error term in the twins conjectures (\ref{twins-a}) or (\ref{twins-b}) is not
known even heuristically).  The only way to
obtain some information about the behaviour of $error ~ term(x,d)$ is to compare
these conjectures with actual computer counts of $\tau_d(x)$. Of course, the best
accuracy has the formula (\ref{main}).  We have compared it with generated by the computer
actual values of $\tau_d(x)$ --- i.e. we have looked at  values of
\bee
\Delta(x,d)\equiv \tau_d(x) - C_2 {\mathfrak S}(d) ~\frac{\pi^2(x)}{x}\Bigg(1-\frac{2\pi(x)}{x}\Bigg)^{{\frac{\mathlarger d}{\mathlarger 2}-1}}.
\eee
The values of $\Delta(x,d)$ were stored for 105 values of $d=2, 4, \ldots,
210(=2\cdot 3\cdot 5\cdot 7)$ at the arguments $x$ forming the geometrical progression
$x_k=1000\times(1.03)^k$. Additionally the values of $|\Delta(x,d)|<9$ were stored
to catch sign changes of $\Delta(x,d)$.
It is difficult to present these data for all values of $d$.  We have found that
for some gaps $d$ there was monotonic
increase of $\Delta(x, d)$, for other gaps there were sign changes of
the difference $\Delta(x, d)$, see Fig.\ref{fig-wykres-delty}. For 30 values of
$d$ of all 105 looked for we have found sign changes for $x<8\times 10^{13}$.
Surprising is the steep growth of $\Delta(x,d)$ for $d=44, 56, 78$ (the same behaviour we
have seen for other values of $d$) in the region of crossing the $y=0$ line.
In fact, there were 76 sign changes of $\Delta(x, 54)$, 109 sign changes of
$\Delta(x, 56)$  and 207 sign changes of $\Delta(x, 78)$.
The general rule is that the ratio
$\tau_d(x)/C_2 {\mathfrak S}(d) ~\frac{\pi^2(x)}{x}(1-\frac{2\pi(x)}{x})^{{\frac{\mathlarger d}{\mathlarger 2}-1}}$
tends to 1, see Fig. \ref{fig-compare_tau_errors}.  Thus we formulate the\\

{\bf Conjecture $\bf 3$}

For every $d$ there are infinitely many sign changes of the functions $\Delta(x,d)$.
For fixed $d$ we guess
\bee
\lim_{x \to \infty} \frac{{\rm Conjecture}_{1,1',1''}(d,x)}{\tau_d(x)} = 1.
\eee

We can test the conjecture (\ref{main2}) with available computer data plotting on
one graph the scaled quantities:
\bee
T_d(x)=\frac{x\tau_d(x)}{C_2 {\mathfrak S}(d) \pi^2(x)}, ~~~~~~~~~~~~~D(x, d)= \frac{d\pi(x)}{x}.
\eee
From   the conjecture (\ref{main2}) we expect that the points
$(D(x, d), T_d(x)), ~~d=2, 4, \ldots, G(x)$ should coincide for each $x$ --- the function
$\tau_d(x)$ displays scaling in the physical terminology. In Fig. \ref{fig-scaling}
we have plotted the points $(D(x, d), T_d(x))$ for $x=2^{28}, 2^{38}, 2^{48}$.
If we denote $u=D(x, d)$ then all these scaled functions should lie on the pure exponential
decrease $e^{-u}$ (Poisson distribution, see \cite[p.60] {Soundararajan-2005}), shown
in red in Fig. \ref{fig-scaling}.
We have determined by the least square method  slopes $s(x)$ of the fits
$a(x)e^{-s(x)u}$ to the linear parts of $(D(x, d), \log(T_d(x)))$. The results are
presented in Fig. \ref{fig-slopes}. The slopes very slowly tend to 1: for over 6
orders of $x$ they change from $1.187$ to $1.136$.

\begin{figure}[h]
\begin{minipage}[t]{0.5\textwidth}
\centering
\includegraphics[width=0.95\textwidth]{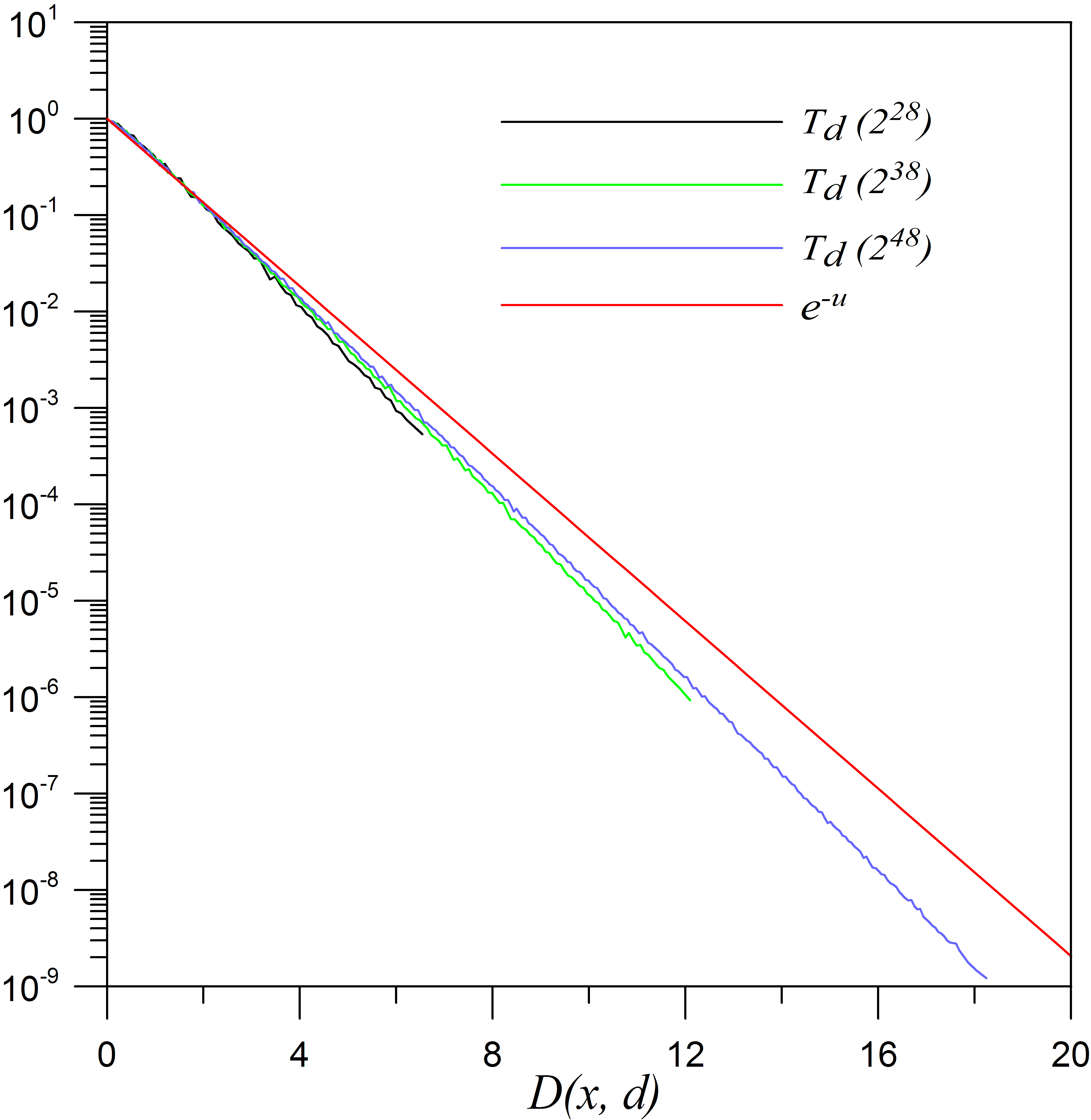}
\caption{Plots of $(D(x, d), T_d(x))$ for $x=2^{28},2^{38},2^{48}$ and in red
the plot of $e^{-u}$. Only the points with $\tau_d(x)>1000$ were plotted to avoid
fluctuations at large $D(x, d)$.}
\label{fig-scaling}
\end{minipage}
\hspace{0.03\textwidth}
\begin{minipage}[t]{0.5\textwidth}
\centering
\includegraphics[width=\textwidth]{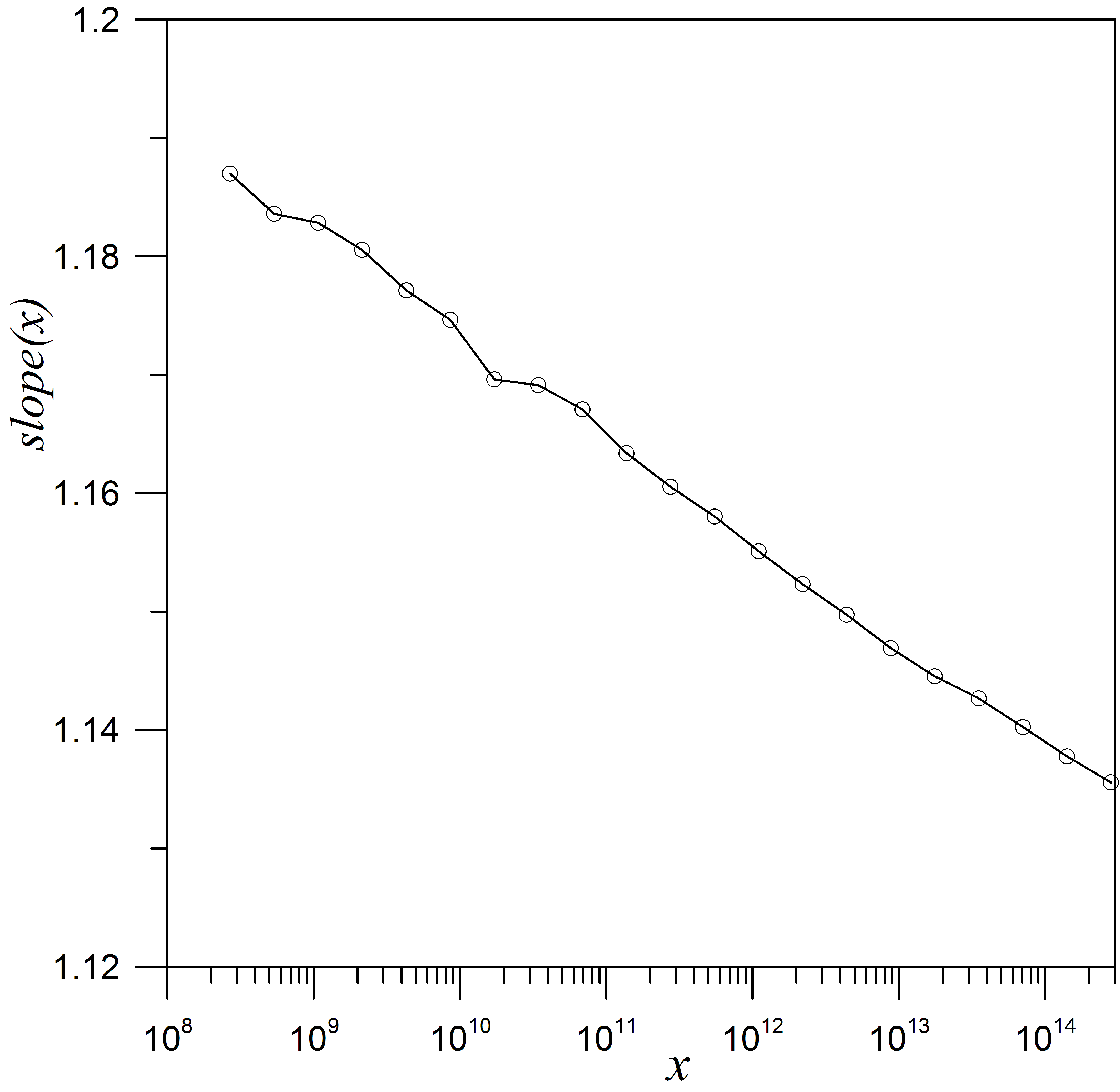}
\caption{ Plot of slopes obtained from fitting straight lines to $(D(d, x), \log(T_d(x)))$
for $x=2^{28},2^{29}, \ldots, 2^{48}$.}
\label{fig-slopes}
\end{minipage}
\end{figure}

\section{Champions -- the most often occurring gaps}
\label{champions}

Direct inspection shows that starting with prime 389 the most often occurring gap between consecutive primes is $d=6$,
with the exception of interval between 491 and 541, where $d=4$  wins over $d=6$, see \cite[Table I]{champions}.
In general the most common gap   between consecutive primes is called champion.
We see in Figure \ref{fig-hist} that for each $x$  the highest value of $\tau_d(x)$  is at  $d=6$. Next local maxima  (spikes)
are at  $d=30$ and $d=210$.  Because the slope  of plots of $\tau_d(x)$ decreases with $x$ at some value $x^{(3)}$ the gap
$d=30=2\cdot 3\cdot 5$ will take  over gap six and at much larger $x^{(4)}$  $d=210=7\sharp$  (product of  first  4  consecutive
primes)   will take over 30. These spikes are caused by the local maxima of the product
${\mathfrak S}(d)$:  when $d$ has many prime factors the product is large.  In fact, in \cite[Lemma 2.2]{Goldston-Leodan-2011}
it is proved, that if $D^{(n)}=2\cdot 3\ldots p_n \equiv p_n\sharp$ denotes the $n$--th primorial, then
${\mathfrak S}(d)<{\mathfrak S}(D^{(n)})$  for every $d$  such that $2<d<D^{(n)}$  (in fact this inequality is fulfilled in the
interval $2<d<D^{(n+1)}$).   From this observation we  can obtain
general estimation for the values of $x^{(n)}$ at which the consecutive products $D^{(n)}= 2\cdot 3\cdot\ldots\cdot p_n$
become the champions. Namely, the condition for appearance of the $n$-th champion is the
following   inequality:
\bee
\tau_{D^{(n-1)}}((x^{(n)})<\tau_{D^{(n)}}((x^{(n)})
\label{warunek}
\eee
From it and (\ref{main2})  ignoring (probably  large)  error terms the estimation  follows:
\bee
\log(x^{(n)}) \sim \frac{2\cdot 3\cdot\ldots \cdot p_{n-1}(p_n-1)}{\log((p_n-1)/(p_n-2))}.
\label{n_ty_champion}
\eee
The values of  turning points for champions obtained from the above formula are presented in Table I.
Obtained from \eqref{n_ty_champion} for $D^{(2)}=6$ the value of
$321$ quite well agrees with the actual value $x^{(2)}=389$.

We will  present the asymptotic form  of \eqref{n_ty_champion} for  large $n$.  We  need  the  closed  formula  for
$2\cdot 3\cdot\ldots \cdot p_{n-1}\cdot p_n$  for  large $n$.   The prime number theorem is sometimes formulated as the law
governing the growth of the Chebyshev function  \cite[Th. 420]{H-W}:
\bee
\vartheta (x)=\sum _{p\leq x}\log(p) \sim x.
\label{Czebyszew}
\eee
We  need  the  value of the  Chebyshev function  at  $n$-th  prime $p_n$  and in \cite[p. 5]{Dusart-2010}  we  find
\bee
\begin{split}
\theta(p_n)=n\log( n) \Bigg( 1+\frac{\log\log( n)}{\log(n)}-\frac{1}{\log(n)} + \frac{\log\log( n)-2}{(\log(n))^2}{\color{white}{\Bigg)}}   \\
{\color{white}{\Bigg(}} -\frac{(\log\log(n))^2- 6\log\log(n)+11}{2(\log(n))^3}+{\mathcal{O}}\left(\frac{(\log\log(n))^3}{(\log( n))^4}\right){\Bigg)}
\end{split}
\eee
Let  us  remark  that   terms  in  the  big  parenthesis are  of  opposite  signs so to   great  extend  they    cancel  out
and    we   will  keep  only  the  first term 1 as  all following  sequences  go to zero with $n \to \infty$:
\bee
\theta(p_n)=n\log( n).   
\eee
So we have the following  rough estimation for primorials (in  fact the error term  is  exponentiated)
\bee
2\cdot 3\cdot\ldots \cdot p_n = p_n\sharp \sim n^n,
\eee
what is an analog of the Stirling  formula for factorials $n!\sim \sqrt{2\pi}n^n e^{-n}$.
For large $n$  we have  $p_n-1\approx p_n,  p_n-2\approx p_n$,   $\log((p_n-1)/(p_n-2))\approx 1/(p_n-2)$ and we obtain
from  \eqref{n_ty_champion}  that $x^{(n)}$  grows asymptotically as:
\bee
x^{(n)} \sim n^{n^{n+1}}.
\eee

\vskip 0.4cm
\begin{table}[ht]
\begin{center}
\bigskip
\begin{tabular}{| c | c | c | c |} \hline
$n$  & $D^{(n)}=p_n \sharp $ & $x^{(n)}$ & $ \prod_{p \mid D^{(n)}, p > 2}\left( {p - 1 \over p - 2}\right) $ \\ \hline
2  & 6 & 3.21$\times 10^{2}$  &     2.00\ldots  \\ \hline
3  & 30 & 1.70$\times 10^{36}$  &     2.67\ldots  \\ \hline
4  & 210 & 5.81$\times 10^{428}$  &     3.20\ldots  \\ \hline
5  & 2310 & 1.48$\times 10^{8656}$  &     3.56\ldots  \\ \hline
6  & 30030 & 1.30$\times 10^{138357}$  &     3.88\ldots  \\ \hline
7  & 510510 & 8.02$\times 10^{3233259}$  &     4.14\ldots  \\ \hline
8  & 9699690 & 8.50$\times 10^{69820169}$  &     4.38\ldots  \\ \hline
9  & 223092870 & 5.14$\times 10^{1992163572}$  &     4.59\ldots  \\ \hline
10  & 6469693230 & 3.56$\times 10^{74595540317}$  &     4.76\ldots  \\ \hline
11  & 200560490130 & 2.10$\times 10^{2486392448589}$  &     4.92\ldots  \\ \hline
12  & 7420738134810 & 1.31$\times 10^{111309396336960}$  &     5.06\ldots  \\ \hline
13  & 304250263527210 & 6.46$\times 10^{5091729308630201}$  &     5.19\ldots  \\ \hline
14  & 13082761331670030 & 1.65$\times 10^{230298784738628635}$  &     5.32\ldots  \\ \hline
\vdots &  \vdots & \vdots & \vdots \\   \hline
\end{tabular} \\
\end{center}
\caption{ The  values of $x^{(n)}$ at which  $D^{(n)}$  become champions (third  column)
and values of the product (\ref{product}) (fourth  column).  }
\end{table}
\vskip 0.4cm

\section{Lemke Oliver-Soundararajan bias}
\label{bias}

In the spring of 2016 some sensation was sparked by the paper \cite{Lemke-Oliver2016}. The authors described biases in the
distribution of pairs of consecutive primes. We will  present their discovery in the particular case of usual base--10  numeral
system, although their  consideration are general.  In the base 10 each prime number (except 2 and 5)   has the last digit
1, 3, 7 or 9,  otherwise it would be  divisible by 2 or 5.  So we have 4 possibilities and  the famous Dirichlet's theorem on
the primes in  arithmetical progressions with the de la Vall{\'e}e Poussin  quantitative supplement, see \cite[chap. 4.V]{Ribenboim},
asserts that all these four possibilities should  be equally probable and one would expect 0.25 \% of primes to end with 1, 3, 7
or 9.  However final digits of the primes that immediately follow them are not equally distributed:  Lemke Oliver
and  Soundararajan found  huge correlation between the last digits of consecutive primes.  For example, we have checked that up
to $2^{34}=17179869184$  the primes with last digit 1 are followed by 36131238 primes with last digit 1, 55962283 primes with last
digit 3,  56247252  primes with last digit 7 and 42391953 primes with last digit 9. The same non--uniform  behavior we obtained for
primes with last digits 3, 7 and 9, see Table 2. We have used there the notation from \cite{Lemke-Oliver2016}: $\pi(x; 10, (a,b)):=
\sharp\{ p_n\leq x:  p_n\equiv a \pmod{10},  p_{n+1}\equiv b \pmod{10} \} $.    In view of the considerations of last
two Sections we should not be surprised by this  outcome:  after primes ending with 1 the  largest number of next primes
should end with 7  as such  primes will follow prime with last  digit 1 after gap 6, 16, 26, ... and the  most  common gap
between primes in this interval  is 6.  Primes  ending with 3 can follow primes  with  last digit 1  after gap 2,
12($=2\cdot 6$), 22, ... etc.  So all these irregularities  reported in  \cite{Lemke-Oliver2016}  are  encoded in the  behavior
of the product  ${\mathfrak S}(d)$ combined with the  exponential with $d$ decay of $\tau_d(x)$.

\vskip 0.4cm

\begin{table}[ht]
\begin{center}
\begin{tabular}{ c || c c c c} \hline
$ a \backslash b $ & $1$ & $ 3 $ & $7$ & $9 $    \\ \hline
\hline
1  &    36131238  & 55962283   & 56247252 & 42391953 \\ \hline
3  &      46218200 & 35046136  & 53225098 & 56247286 \\ \hline
7  &       48527277 & 51208545  & 35035483 & 559654613 \\ \hline
9  &       59856011 & 48519757  & 46228932  & 36128195  \\   \hline
\end{tabular}
\end{center}
\caption{ The $4\times 4 $ matrix  of  values of $\pi(2^{34}; 10, (a,b))$,  $a,b =1,3,7,9$.  The reader can explain why the sum of
all entries above is by 4 smaller than  $\pi(2^{34})=762939111$. }
\end{table}
\vskip 0.4cm

\section{Maximal gap between consecutive primes}
\label{G-maxy}

From (\ref{main}) or  (\ref{main2}) we can obtain approximate formula for  $G(x)$
assuming that maximal difference $G(x)$  appears only once, so
$\tau_{G(x)}(x)=1$: simply the largest gap is equal to the value of $d$ at
which $\tau_d(x)$ touches the $d$-axis on  Fig.\ref{fig-hist}. Skipping the
oscillating term ${\mathfrak S}(d)$, which is very often close to 1,
we get for $G(x)$ the following estimation expressed directly by $\pi(x)$:

{\bf Conjecture $\bf 4$}
\bee
G(x) \sim g(x)\equiv \frac{x}{\pi(x)} \big(2\log( \pi(x)) -\log(x) + c\big),
\label{d_max}
\eee
where $c=\log(C_2)= 0.2778769\ldots$.

The above formula explicitly reveals the fact that the value
of $G(x)$ is connected with the number of primes $\pi(x)$: more primes means smaller
$G(x)$. For the Gauss approximation $\pi(x)\sim x/\log(x)$  the following dependence follows:
\bee
G(x) \sim \log(x) (\log(x)-2\log\log(x) + c)
\label{d_max2}
\eee
and for large $x$ it passes into the well known  Cramer's   conjecture \cite{Cramer}:
\bee
G(x)\sim \log^2(x).
\label{eCramer}
\eee
The examination of the formula (\ref{d_max})  and the formula (\ref{eCramer})
with the available results of the computer search  is given in Fig.\ref{fig-G}.
The lists of known maximal gaps between consecutive primes we have
taken from our own computer search up to $2^{48}$ and larger from web sites
www.trnicely.net  and  www.ieeta.pt/$\sim$tos/gaps.html.
The largest known  gap 1510  between consecutive primes follows the prime
$6787988999657777797=6.788 \ldots\times10^{18}$. On these web sites tabulated  values of
$\pi(x)$ can  also  be  found and we have used them to plot the formula (\ref{d_max}).
Let $\nu_G(T)$ denotes the number of sign changes of the difference $G(x)-g(x)$
for $2<x<T$. There are 33 sign changes of the difference $G(x)-g(x)$
in the  Fig.\ref{fig-G}  and $\nu_G(T)$ is  presented in the inset in Fig. \ref{fig-G}.
The least  square method gives for $\nu_G(T)$  the equation $0.78\log(T) + 0.63$.

There appeared in literature a few other formulas for approximate values of  $G(x)$, see
e.g.  \cite{Shanks1964}, \cite{Cadwell}; in particular
D.R. Heath-Brown in \cite[p. 74]{Heath-Brown-1988} gives the following formula:
\bee
G(x) \sim \log(x)(\log(x)+\log\log\log(x)).
\eee

A. Granville argued  \cite{Granville} that the actual $G(x)$ can be larger
than that given by Cramer's  model  (\ref{eCramer}), namely he claims  that
there are infinitely many pairs of primes $p_n, p_{n+1}$ for which:
\bee
p_{n+1} - p_n = G(p_n)> 2 e^{-\gamma}\log^2(p_n) = 1.12292\ldots \log^2(p_n).
\label{Granville-Cramer}
\eee
where $\gamma=0.577216\ldots$ is the Euler--Mascheroni constant.
The estimation (\ref{Granville-Cramer}) follows from the inequalities proved by
H.Maier in the paper \cite{Maier-1985}, which put into doubts  Cramer's ideas.
For other contradiction between Cramer's model and the  more strict results, see \cite{Pintz-2007}.

\begin{figure}[ht]
\begin{center}
\includegraphics[width=\textwidth,angle=0]{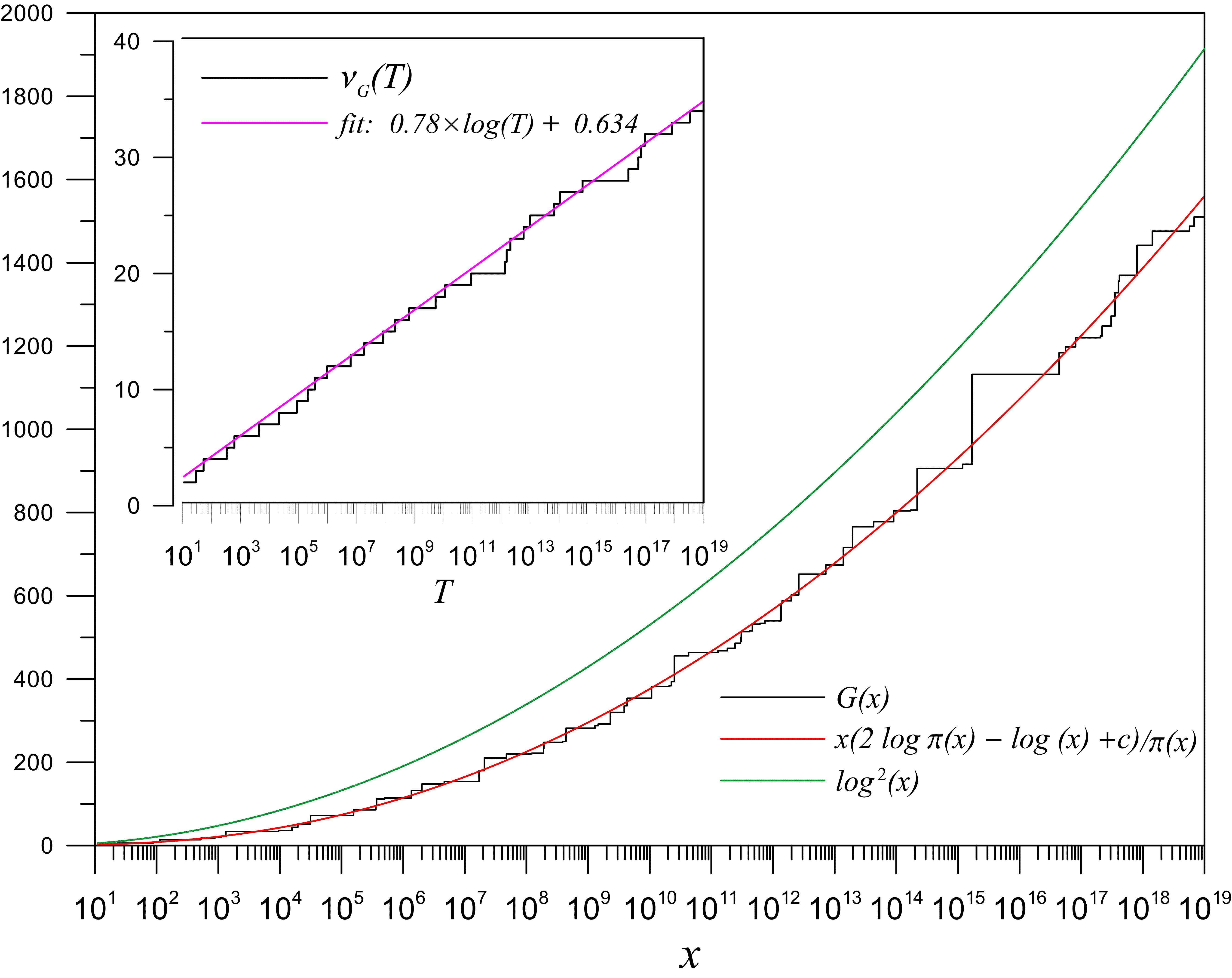}
\vspace{0.2cm}
\caption{The comparison of $G(x)$ and $g(x)$ as well as of
the Cramer conjecture $G(x) \sim \log^2(x)$. In the inset there is the plot of a number of crossings of curves
representing $G(x)$ and $g(x)$.  This figure  should be compared with the figure on
page 12 in \cite{Zagier1977}. }
\label{fig-G}
\end{center}
\end{figure}

\section{Generalized Brun's constants}
\label{Bruny}

In  1919 Brun \cite{Brun} has shown that the sum of the  reciprocals of all twin primes is finite:
\bee
\BB_2=\left( \frac{1}{3} + \frac{1}{5}\right) + \left(\frac{1}{5} + \frac{1}{7}\right) + \left(\frac{1}{11} +
\frac{1}{13}\right) + \ldots < \infty.
\label{def_B2}
\eee

Sometimes 5 is included only once, but here we will adopt the above  convention. The analytical formula for $\BB_2$ is unknown
and the value of the sum (\ref{def_B2}) is called the Brun's constant \cite{Wrench1974}. The numerical estimations give
\cite{Nicely_Brun} $\BB_2=1.90216 058\ldots$.   Here we are going to generalize the above $\BB_2$ to the sums
of reciprocals of all consecutive primes separated by gap $d$ and to  propose a compact expression giving the values of
these sums for $d\geq 6$.\\

Let ${\cal T}_d$ denote the set of consecutive primes
separated by distance $d$:
\bee
{\cal T}_d = \{(p_{n+1}, p_n) : p_{n+1} - p_n=d\}.
\eee
We define the  generalized Brun's constants by the formula:
\bee
{\cal B}_d=\sum_{p \in {\cal T}_d}  \frac{1}{p}.
\label{def_Bd}
\eee
We adopt the rule, that if a given gap $d$ appears two times in a row:
$p_n-p_{n-1}=p_{n+1}-p_n$,  the corresponding
middle prime $p_n$  is counted two times (in the case of $\BB_2$ only 5 appears two times);
e.g. for $d=6$ we have the terms $\ldots + 1/47+ 1/53+ 1/53 +1/59+\ldots$ and next
$\ldots + 1/151+1/157+1/157+1/163+\ldots$.

B.Segal has proved \cite{Segal} that the sum in (\ref{def_Bd}) is convergent for
every $d$, thus generalized Brun's constants are finite. Because of that the sums
(\ref{def_Bd}) can be  called Brun--Segal constants for $d>2$.\\

Let us define partial (finite) sums:
\bee
{\cal B}_d(x)=\sum_{p \in {\cal T}_d, p<x} \frac{1}{p}.
\eee

We have computed on  the computer quantities $\BB_d(x)$ for $x$ up to
$x=2^{46}\approx 7.037\times 10^{13}$.  
In  Fig. \ref{sample-Bruny} we present a part of the obtained data.

\begin{figure}[ht]
\begin{center}
\includegraphics[width=\textwidth,angle=0]{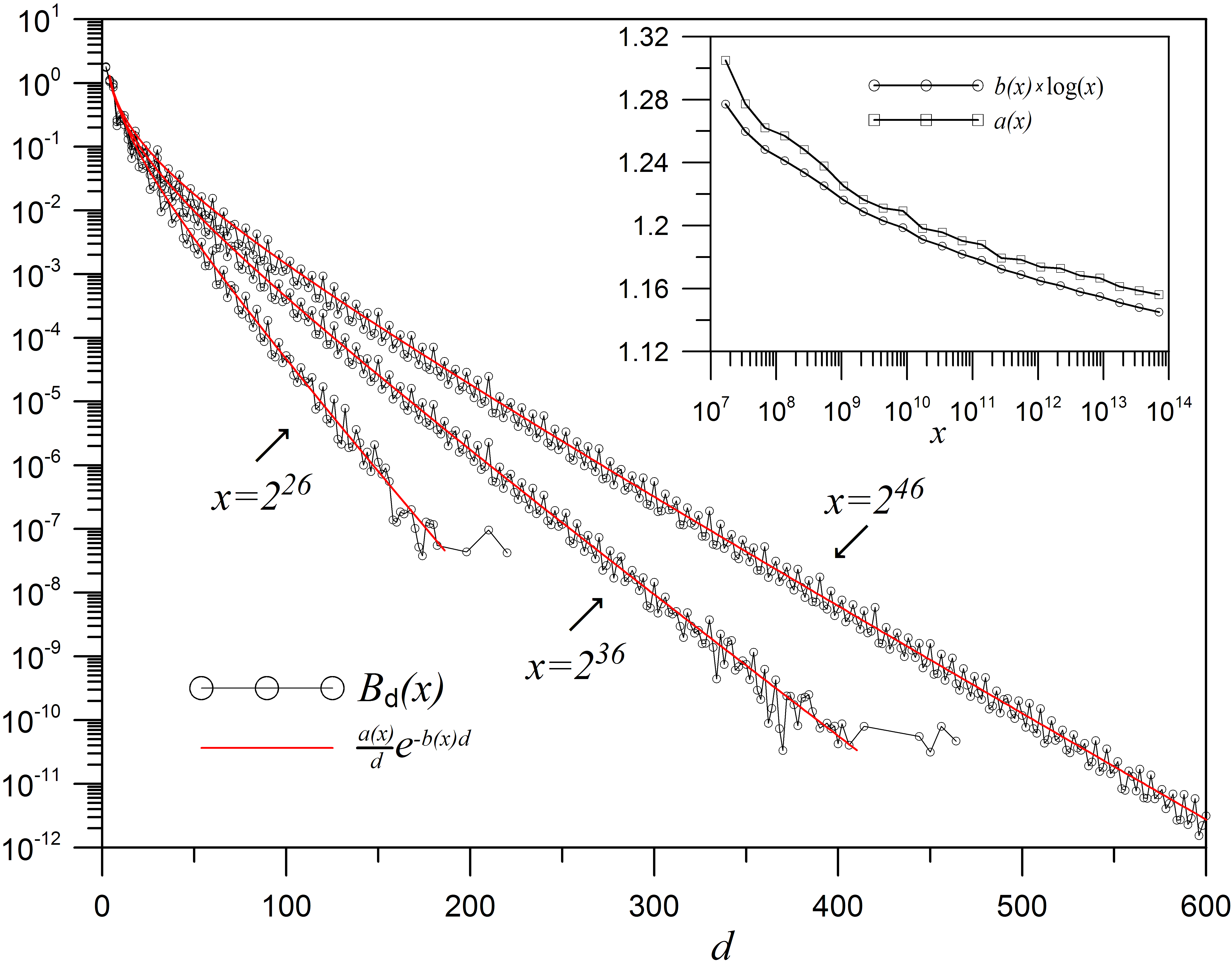} \\
\caption{ The plot of $\BB_d(x)$ for $x=2^{26}=6.71...\times 10^7,
2^{36}=6.87...\times 10^{10}, 2^{46}=7.04...\times 10^{13}$. The fit $a(x)e^{-b(x)d}/d$
to $\BB_d(x)/{\mathfrak S}(d)$ obtained by the least square method is plotted in red.
In the inset the values of $a(x)$ as well as the ratio  between conjectured slope
$-1/\log(x)$ and actual fit $b(x)$ are shown for $x=2^{24}, 2^{25}, \ldots, 2^{46}$.}
\label{sample-Bruny}
\end{center}
\end{figure}

The dependence of $\BB_2(x)$ on $x$ is usually (see \cite{Wrench1974},
\cite{Brent}) obtained by appealing to the conjecture (\ref{twins-a})
(i.e. Hardy--Littlewood conjecture (\ref{H-L}) for $d=2$).
It gives that the probability
to find a pair of twins in the vicinity of $x$ is $2c_2/\log^2(x)$, so
the expected value of the finite approximation to the Brun constant can be
estimated as follows:
\bee
\BB_2(x) = \BB_2(\infty) - \sum_{p\in \TT_2, p>x} \frac{1}{p} \approx
\BB_2 - 2C_2\int_x^\infty  \frac{du}{u\log^2(u)} = \BB_2 - \frac{2C_2}{\log(x)}.
\label{B2}
\eee
It means that the plot of finite approximations $\BB_2(x)$ to the  original Brun constant is a linear function of
$1/\log(x)$  and intercept for $x=\infty$ of this  plot of $\BB_2(x)$  vs $1/\log(x)$  gives $\BB_2$. In other words,
the value of  $\BB_2$  is extrapolated from finite  sum $\BB_2(x)$  by adding to it term  $2C_2/\log(x)$. The same reasoning
applies  {\it mutatis mutandis} to the gap $d=4$.

To repeat the above reasoning for $d=2, ~4$ for larger $d$ an analog
of the Hardy--Littlewood conjecture for the pairs of {\it consecutive} primes
separated by distance $d$ is needed and we will use the form (\ref{main2}) for $\tau_d(x)$
(the integrals occurring below can be calculated analytically also for (\ref{main})).
Putting in the equation (\ref{main2}) $\pi(x)=x/\log(x)$ we obtain for $d\geq 6$ :
\bee
\BB_d(x) = \BB_d(\infty) - \sum_{p\in \TT_d, p>x} \frac{1}{p} \approx
\BB_d - 2C_2  \prod_{p\mid d} \frac{p-1}{p-2}\int_x^\infty \frac{e^{-d/\log(u)}}{u\log^2(u)} du.
\eee
and the integral can be calculated explicitly:
\bee
\BB_d(x) \approx  \BB_d(\infty) + \frac{2C_2}{d} \prod_{p\mid d} \frac{p - 1}{p - 2} \left(e^{-d/\log(x)}-1\right).
\label{Bd_x}
\eee
From this, it follows that the partial sums  $\BB_d(x)$ for $d\geq 6$
should depend linearly on $e^{-d/\log(x)}$ instead of linear dependence on
$1/\log(x)$ for $\BB_2(x)$ and $\BB_4(x)$.\\

Because $\BB_d(x)$ is 0 for $x=1$ (in fact each $\BB_d(x)$ will be zero up to
the first occurrence of the gap $d$), we take in (\ref{Bd_x})  the limit $x\to 1^+$ and obtain

{\bf Conjecture $\bf 5$}
\bee
\BB_d(\infty) \equiv \BB_d \approx \frac{2C_2}{d} \prod_{p\mid d} \frac{p - 1}{p - 2} ~~{\rm for~ }d\geq 6.
\label{main-B-d}
\eee
Thus the formula expressing the $x$ dependence of $\BB_d(x)$ has the form:
\bee
\BB_d(x) \sim  \frac{2C_2}{d}  \prod_{p\mid d} \frac{p - 1}{p - 2} e^{-d/\log(x)}.
\label{Bd-od-x}
\eee

The characteristic shape of the dependence of  $\BB_d(x)/{\mathfrak S}(d)$  on $d$
is described by the relation $\log(\BB_d(x)/{\mathfrak S}(d))\sim -\log(d) - d/\log(x)$:
if $d/\log(x)>\log(d)$ the linear dependence on $d$ preponderates. We have fitted by least square method the
dependence $\log(a(x))-db(x)$ to the actual values of $\log(d \BB_d(x)/2C_2{\mathfrak S}(d))$. We  obtained,
that indeed $b(x)$ tends to $1/\log(x)$ and $a(x)$ tends to 1 with increasing $x$, see the  inset in Fig. \ref{sample-Bruny}.

The comparison of the formula (\ref{main-B-d}) with the values extrapolated
from the partial  approximations $\BB_d(2^{46})$
\bee
\BB_d(\infty) = \BB_d(2^{46}) + \frac{2C_2}{d} \prod_{p\mid d} \frac{p - 1}{p - 2}
\left(1-e^{-d/46\log(2)}\right)
\label{extrapolated}
\eee
obtained from the equation (\ref{Bd_x}), is shown in  Fig. \ref{Bd-extrapolated-and-ratio}
for $d\geq 6$ --- predicted by (\ref{main-B-d}) values for $d=2$ and $d=4$ are skipped.
Because on average the product ${\mathfrak S}(d)$ is equal to $1/c_2$, we can write $\BB_d\approx 4/d$.
Let us mention that  $4/d$ provides remarkably good approximations to
$\BB_2=1.90216058\ldots$ and $\BB_4=1.19705\ldots$.

\begin{figure}[ht]
\begin{center}
\includegraphics[width=\textwidth,angle=0]{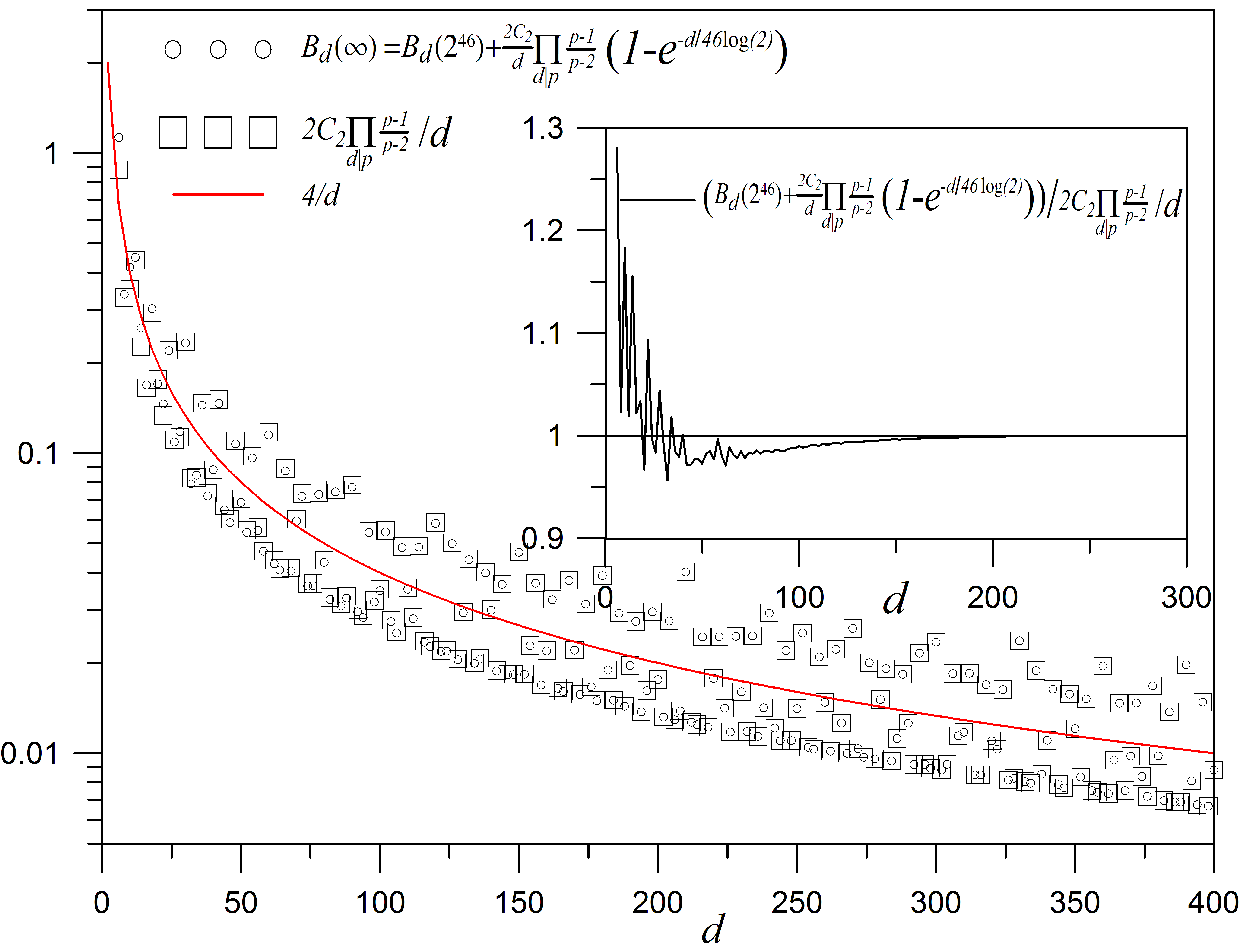} \\
\caption{The plot of the generalized Brun's constants $\BB_d$ extrapolated from
(\ref{extrapolated}) (marked by circles)  and predicted  by (\ref{main-B-d}) (marked by
squares). In the inset the ratio of the values obtained from these two equations is
plotted.}
\label{Bd-extrapolated-and-ratio}
\end{center}
\end{figure}

\section{The Merten's Theorem on the prime harmonic sum.}
\label{Mertens}

Leonhard  Euler showed that the sum  of reciprocals  of all primes  $p<x$  diverges  like $\log(\log(x))$  and it was the  first
constructive proof of infinitude of primes.  In 1874 F. Mertens proved  more precise dependence
 \cite{Mertens-1874}, \cite[Theorems 427 and 428]{H-W}, \cite{Villarino-Mertens}:
\bee
\sum_{p<x} \frac{1}{p} = \log(\log(x)) + M + o(1);
\label{suma_p}
\eee
here $M=0. 2614972\ldots$ is the Mertens constant which has a few representations:
\bee
M=\sum_p (\log(1-1/p)+1/p)= \gamma + \sum_{k=2}^\infty \mu(k)\log(\zeta(k))/k,
\label{eq-Mertens}
\eee
where $\mu(n)$ is the Moebius function and $\zeta(s)$ is the Riemann zeta function.  In Fig. \ref{fig-Mertens} we present
comparison of the above formula with data from our computer calculation of $\BB_d(x)$.
On the other hand, the sum $\sum_{p<x} 1/p$ can be expressed by finite approximations
to the generalized Brun's constants:
\bee
\sum_{p<x} \frac{1}{p} = {1\over 2} + {1\over 6} + {1\over 2} \sum_d
\BB_d(x) = M^\prime + {2\over 3} + C_2 \sum_{d=2}^{G(x)} \frac{1}{d} \prod_{p\mid d} \frac{p - 1}{p - 2} e^{-d/\log(x)}
\label{sum1}
\eee
Because each prime except 2 and 3 (hence the terms $1/2$ and $1/6$) appears
as the right and left end of the adjacent pairs, we have to divide the sum by  2
(we remind that we have adopted in previous Section  the convention that if
a given gap $d$ appears two times in a row: $p_n-p_{n-1}=p_{n+1}-p_n$ --- the
corresponding  middle prime $p_n$  is counted two times).
We have introduced above the constant $M^\prime$ which accounts the sum of
the unknown errors terms in (\ref{Bd-od-x}) as well as incorporates the
fact that  the dependence of $B_2(x)$ and $B_4(x)$ on $x$ is not described by
the formula (\ref{Bd-od-x})   but by (\ref{B2}).
The sum in (\ref{sum1}) runs over even $d$ and extends up to the greatest
gap $G(x)$ between two consecutive primes smaller than $x$. For $G(x)$ we will
use the Cramer's formula (\ref{eCramer}):  $G(x)\approx \log^2(x)$.
To get rid of the product ${\mathfrak S}(d)$, we will make use of the (\ref{rownosc1})
and we obtain:
\bee
\sum_{p<x} \frac{1}{p} = M^\prime + {2\over 3} + 2 \sum_{d=2}^{G(x)} \frac{1}{d}e^{-d/\log(x)} =  M^\prime +  \frac{2}{3} +
\sum_{k=1}^{\frac{1}{2}G(x)} \frac{1}{k} q^k,~~~ q=e^{-2/\log(x)}.
\label{sum2}
\eee
Expanding $\log(1-q)$, where $0<q<1$, into the series we obtain
\bee
\sum_{k=1}^n {1\over k} q^k = -\log(1-q) + \int_0^q \frac{u^n}{u-1} du.
\label{tozsamosc}
\eee
For large $x$ the term  with logarithm goes into:
\bee
\log(1 - e^{-2/\log(x)}) = - \log(\log(x)) + \log(2) + \OO(1/\log(x)).
\label{sum4}
\eee
Now, by the weighted  mean value theorem  we calculate the integral:
\bee
\II = \int_0^q \frac{u^n}{u-1} du = \frac{1}{(\theta q - 1)} \frac{q^{n+1}}{(n+1)},~~~0<\theta<1.
\eee
But $q=\exp(-2/\log(x))<1$ and:
\bee
\Big| \frac{1}{\theta q - 1}\Big| < \frac{1}{1 - q} = \frac{e^{2/\log(x)}} {e^{2/\log(x)}-1} < \frac{\log(x)}{2} e^{2/\log(x)} = \OO(\log(x)).
\eee
For large  $x$ we have in the virtue of the Cramer conjecture that in our
case $n\sim {1\over 2}\log^2(x)$, thus  we  have  heuristically:
\bee
\mid \II \mid = \OO (1/x\log(x)).
\eee

\begin{figure}[ht]
\vspace{-1.1cm}
\hspace{-2.5cm}
\begin{center}
\includegraphics[width=\textwidth, angle=0]{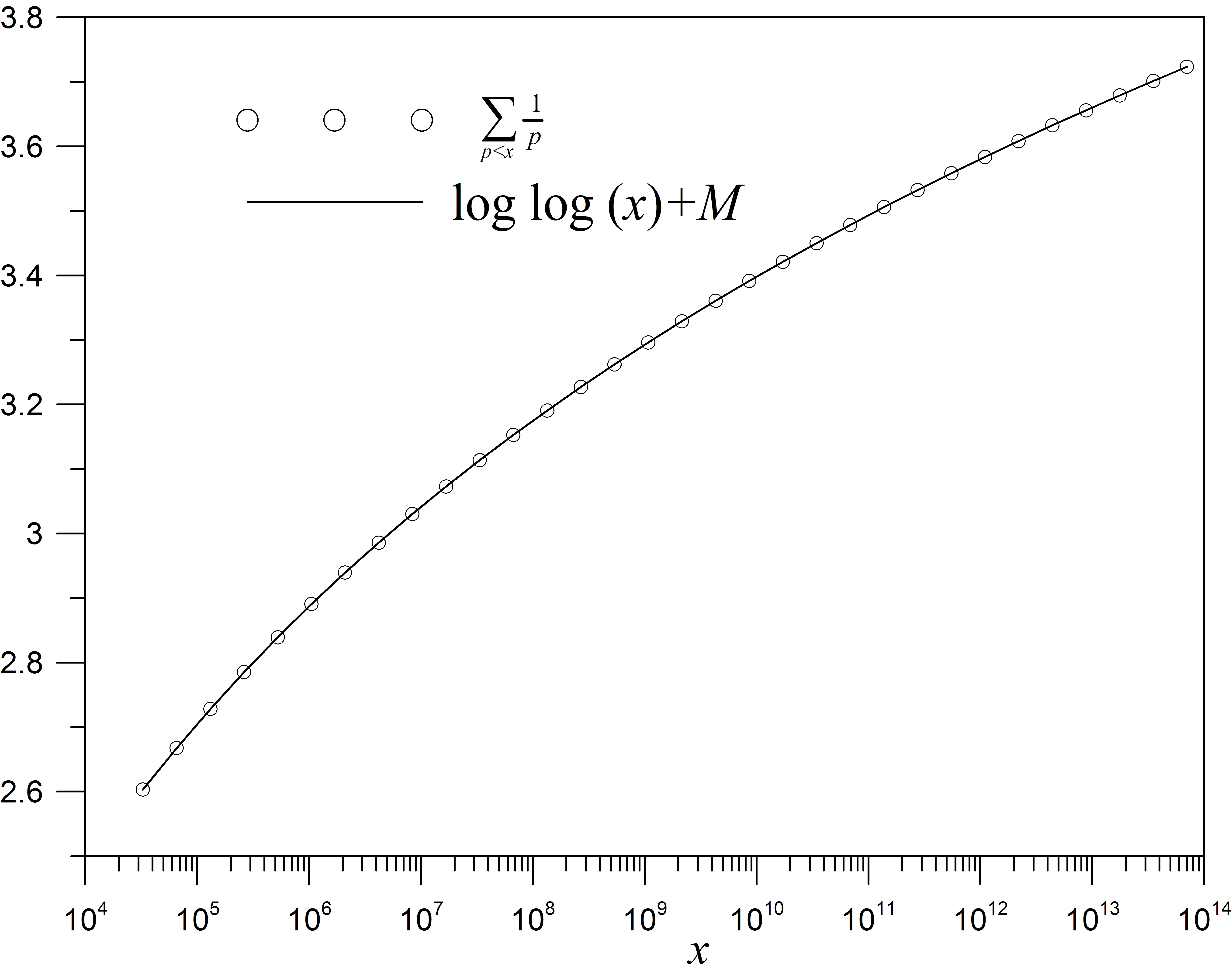} \\
\vspace{0.0cm}
\caption{\small The plot of the prime harmonic sum up to $x=2^{15}, 2^{16}, \ldots,
2^{46}$ and the Merten's approximation  to it. The original of this figure has
$y$ axis of the length 8 cm and spans the interval (2.5, 3.8), so if the $x$ axis would
be plotted in the linear scale instead of logarithmic, then it should be
$5.33(3)\times 10^9$ km long ---  that is the size of  the Solar System.}
\label{fig-Mertens}
\end{center}
\end{figure}

Finally we obtain from (\ref{Bd-od-x}) and (\ref{sum1}):
\bee
\sum_{p<x} \frac{1}{p} = \log(\log(x)) + M^\prime + \frac{2}{3}  - \log(2) + \OO(1/\log(x))
\label{sum5}
\eee
Because 2/3 is practically equal to $\log(2)$ to require consistency with
the Merten's theorem,  we have to postulate that $M^\prime \approx M$.

\section{First occurrence of a given gap between consecutive primes}
\label{first-d}

In this section we will present  the heuristical reasoning leading to the formula
for the first appearance of a given gap of length $d$,  see e.g. \cite{Lander-Parkin},
\cite{Brent1980}, \cite{Young-Potler}, \cite{Nicely1999}.

We will use the conjecture (\ref{Bd-od-x}) to estimate
the position of the first appearance of a pair of primes separated by a
gap of the length $d$. More specifically, let:
\bee
p_f(d)=\begin{cases} {\rm minimal~prime, ~such~ that~the~ next~ prime~}
p^\prime=p_f(d) + d \\
 \infty \text{ if there is no pair of primes } p_{n+1}-p_n=d. \\
\end{cases}
\eee
It is not known whether gaps of arbitrary  length $d$ exist or
not, i.e.  whether  for every even $d$ there  is $p_f(d) < \infty$   \cite{Brent1980} (consult the Polignac's conjecture).

We can obtain the heuristic formula for $p_f(d)$ by remarking  that the finite
approximations to the generalized Brun's constants are for the first time
different from zero at $p_f(d)$ and then they are equal to $2/p_f(d)$:
\bee
\frac{4c_2}{d} \prod_{p\mid d} \frac{p - 1}{p - 2} e^{-d/\log(p_f(d))}=\frac{2}{p_f(d)}.
\label{eq-p_f}
\eee
Referring to the argument that on average ${\mathfrak S}(d)$ is equal to $1/c_2$, we skip  ${\mathfrak S}(d)$ and $c_2$.   Neglecting the $\log(2)=0.69314\ldots$, we end up with the quadratic equation for $t=\log(p_f(d))$:
$$
t^2 - t\log(d) - d = 0
$$
The positive solution of this equation gives

{\bf Conjecture $\bf 6$}
\bee
p_f(d) \sim \sqrt{d}\h{ \rm e}^{\frac{1}{2}\sqrt{\log^2(d) + 4d}}.
\label{main-p_f}
\eee

The comparison of this formula with the actual available data from the computer search
is shown in Fig. \ref{fig-first-d}. Most of the points plotted on this figure come
from our  own search up to $2^{48}=2.815 \ldots \times 10^{14}$.
First occurrences $p_f(d) > 2^{48}$ we have taken from http://www.trnicely.net
and  http://www.ieeta.pt/$\sim$tos/gaps.html. In the Fig.\ref{fig-first-d} there is
also a plot of the conjecture made by Shanks \cite{Shanks1964}:
\bee
p_f(d)\sim {\rm e}^{\sqrt{d}},
\label{Shanks1}
\eee
while from (\ref{main-p_f}) for large $d$ it follows that
\bee
p_f(d)\sim \sqrt{d} \h{\rm e}^{\sqrt{d}}.
\label{moje-p-f}
\eee

\vspace{0.4cm}

\begin{figure}[ht]
\begin{center}
\includegraphics[width=\textwidth, angle=0]{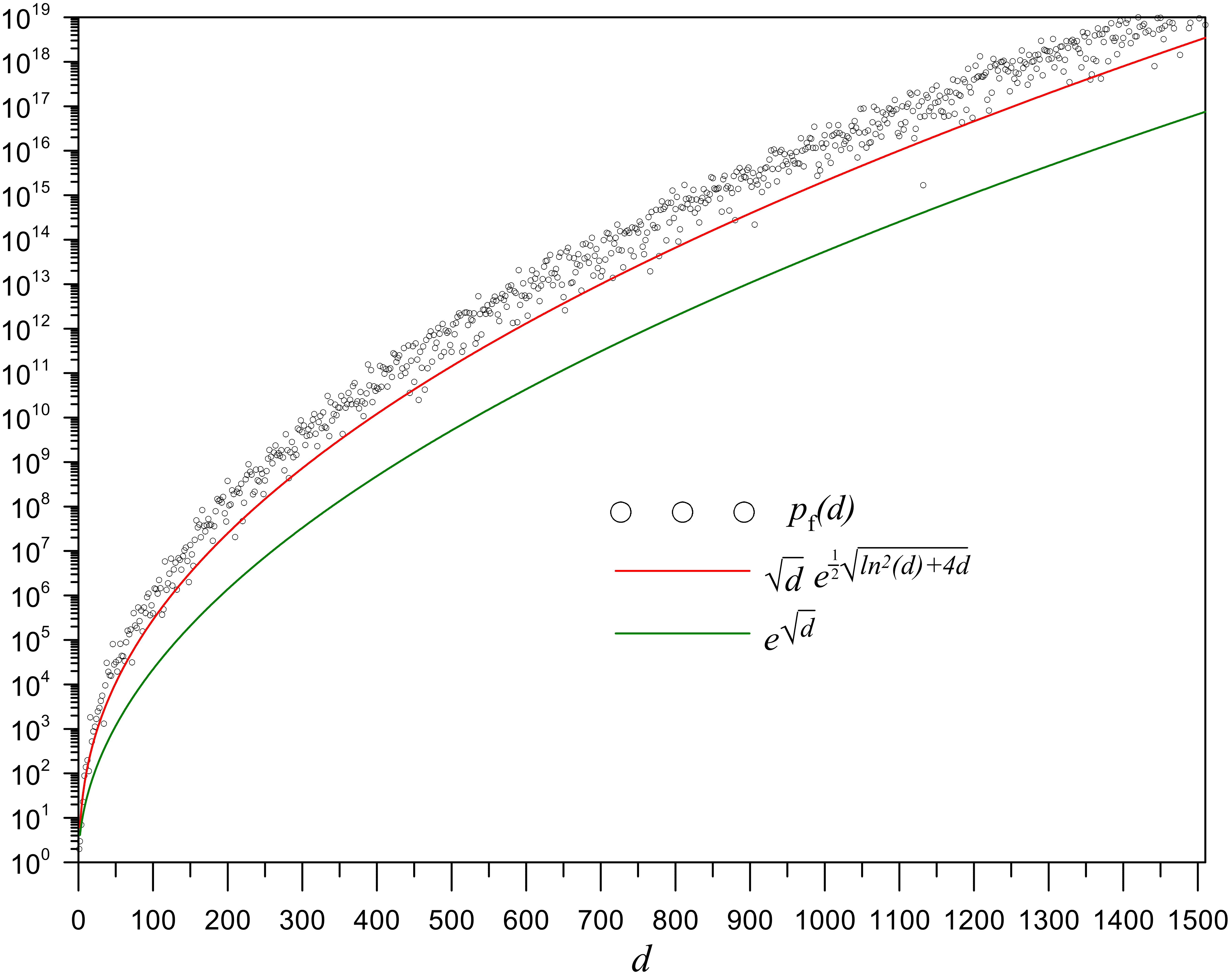}
\vspace{0.0cm}
\caption{ The plot of $p_f(d)$ and approximation to it given by (\ref{main-p_f}) and (\ref{Shanks1}).}
\label{fig-first-d}
\end{center}
\end{figure}

As an application of (\ref{d_max}) we can answer the question
raised in 1964 by P.A.Carlson, who wanted to know at which order of
magnitude $N$ the first appearance of the string of all one million
consecutive numbers being composite can be found. A very rough estimate
was found by D.Shanks \cite{Shanks1964}. Shanks found:
\bee
10^{300} < p_f(10^6) < 100^{600}.
\eee
From  \eqref{moje-p-f}  we obtain  $p_f(10^6)\approx 1.97 \times 10^{437}$.

We make remark concerning the Cramer's formula for maximal gap $G(x)\sim \log^2(x)$. Namely reverting it we obtain that
the maximal gap $g$  appears at $x\sim exp(\sqrt{g})$,  while arbitrary gap $d$ appears for the first time at
$\sqrt{d} \h{\rm e}^{\sqrt{d}}$.  Hence the maximal gaps are those gaps which appear for the first time at $x$ roughly
by $\sqrt d$ earlier than  the  typical  size $p_f(d)\sim \sqrt{d} \h{\rm e}^{\sqrt{d}}$.

\section{The Andrica Conjecture}
\label{Andrica-conj}

In the last section we will make use of most of the conjectures formulated so far.
The Andrica conjecture \cite{Andrica} (see also \cite[p. 21]{Guy} and
\cite[p. 191]{Ribenboim}) states that the inequality:
\bee
A_n \equiv \sqrt{p_{n+1}} - \sqrt{p_n} < 1,
\label{Andrica-ineq}
\eee
holds for all $n$.  Despite its simplicity
it remains unproved.  In Table 3 the values of $A_n$ are
sorted in descending order  (it is believed this order will persist forever).   We have
\bee
\sqrt{p_{n+1}} - \sqrt{p_n}=\frac{p_{n+1} - p_n}{\sqrt{p_{n+1}} + \sqrt{p_n}}<
\frac{d_n}{2\sqrt{p_n}}.
\label{nierownosc}
\eee
From this we see that the growth rate of the form $d_n = \OO(p_n^\theta)$ with
$\theta<1/2$ will suffice for the proof of (\ref{Andrica-ineq}), but as we have
mentioned in the Introduction, currently the best unconditional result is
$\theta=21/40$ \cite{Baker-et-al}.  If \eqref{Andrica-ineq} is true then $d_n=\OO(\sqrt{p_n})$ and the Legendre  conjecture that
between $n^2$ and $(n+1)^2$  there is always  a  prime  follows.  Put in another words: if $A_n <1 $, then there must be
a prime between $n^2$ and $(n+1)^2$.

For twins primes $p_{n+1}=p_n+2$ there is no problem
with (\ref{Andrica-ineq}) and in general for short gaps $d_n=p_{n+1} - p_n$
between consecutive primes the inequality (\ref{Andrica-ineq}) will be satisfied.
The Andrica conjecture can be violated only by  extremely large gaps between
consecutive primes. Let us denote the pair of primes $<x$  comprising the largest
gap $G(x)$ by $p_{L+1}(x)$ and $p_L(x)$, hence we have
\bee
G(x)=p_{L+1}(x)- p_L(x).
\eee
Thus we will concentrate on the values of the difference appearing in
(\ref{Andrica-ineq}) corresponding to the largest gaps and so let us
introduce the function:
\bee
R(x)=\sqrt{p_{L+1}(x)}- \sqrt{p_L(x)}.
\eee
Then we have:
\bee
A_n\leq R(p_n).
\eee

\vskip 0.4cm
\begin{table}[ht]
\begin{center}
\begin{tabular}{|c | c | c | c | c |} \hline
$ n $ & $p_n$ & $ p_{n+1} $ & $ d_n $ & $\sqrt{p_{n+1}}-\sqrt{p_n} $    \\ \hline
4  &       7 & 11  & 4 & 0.6708735 \\ \hline
30  &      113 & 127  & 14 & 0.6392819 \\ \hline
9  &       23 & 29  & 6 & 0.5893333 \\ \hline
6  &       13 & 17  & 4 & 0.5175544 \\ \hline
11  &      31 & 37  & 6 & 0.5149982 \\ \hline
2  &       3 & 5  & 2 & 0.5040172 \\ \hline
8  &       19 & 23  & 4 & 0.4369326 \\ \hline
15  &      47 & 53  & 6 & 0.4244553 \\ \hline
46  &      199 & 211  & 12 & 0.4191031 \\ \hline
34  &      139 & 149  & 10 & 0.4167295 \\ \hline
\vdots &  \vdots & \vdots & \vdots & \vdots \\   \hline
\end{tabular} \\
\end{center}
\caption{The values of $\sqrt{p_{n+1}}-\sqrt{p_n}$  sorted in descending order.}
\end{table}
\vskip 0.4cm

\begin{figure}[ht]
\begin{center}
\includegraphics[width=0.9\textwidth, angle=0]{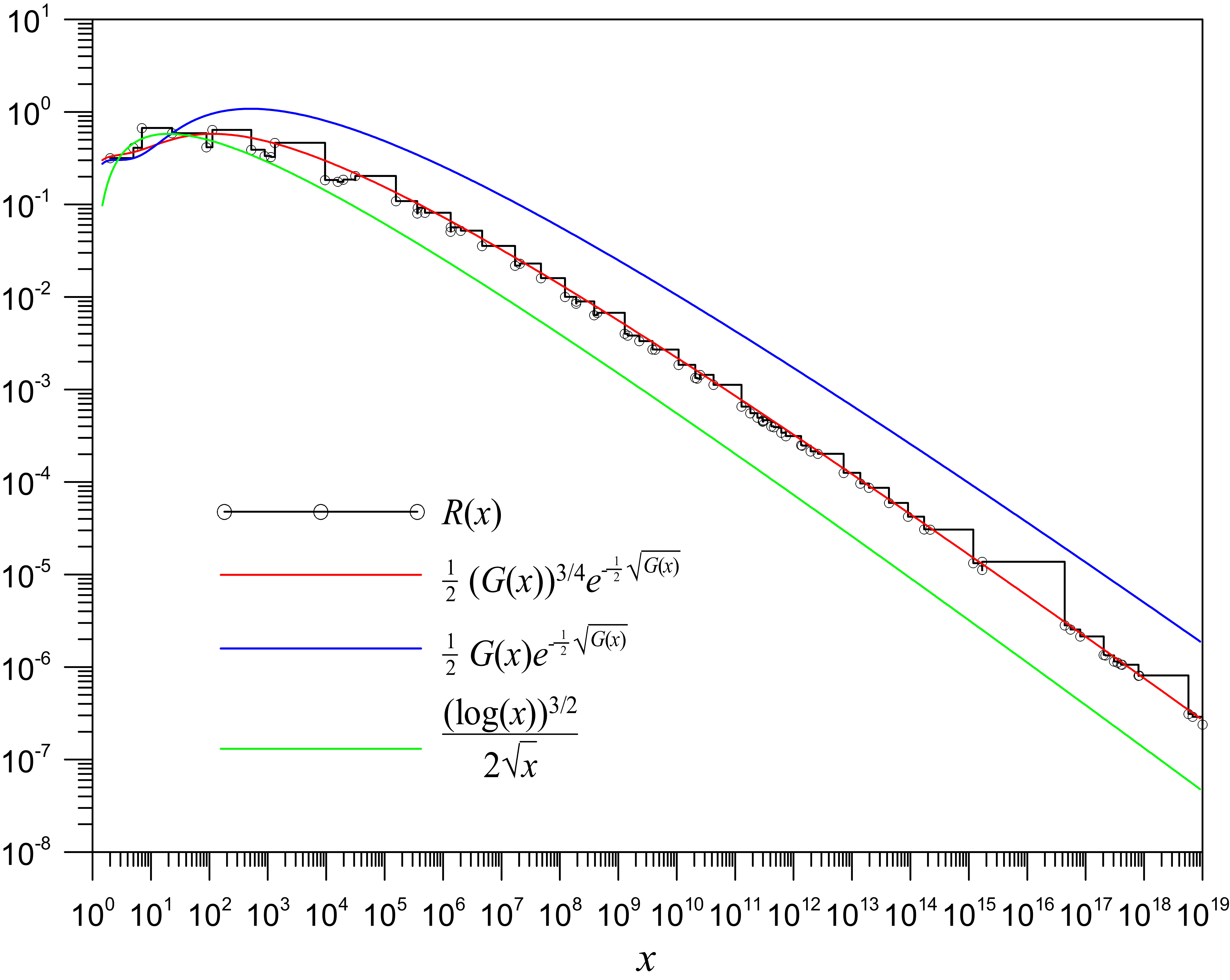} \\
\vspace{0.0cm}
\caption{\small The plot of $R(x)$ and approximations to it
given by (\ref{main-formula}), (\ref{R-Cramer}) and (\ref{Shanks-2}). There are 77 maximal gaps available currently and
hence there are 77 circles in the plot of $R(x)$. To calculate $g(x)$ given by
(\ref{d_max}) we have used tabulated values of $\pi(x)$ available at the web sites
www.trnicely.net  and  www.ieeta.pt/$\sim$tos/primes.html.  There are over 50 crossings  (sign changes of difference)
of our formula (\ref{main-formula}) with $R(x)$.}
\label{fig-Andrica}
\end{center}
\end{figure}

For a given gap $d$ the largest value of the difference $\sqrt{p+d}-\sqrt{p}$
will appear at the first appearance of this gap:
each next pair $(p', p'+d)$ of consecutive primes separated by $d$ will
produce smaller difference (see (\ref{nierownosc})):
\bee
\sqrt{p'+d}-\sqrt{p'}<\sqrt{p+d}-\sqrt{p}.
\eee
Hence, we have to focus our attention on the first occurrences $p_f(d)$ of the gaps.
Using the conjecture (\ref{moje-p-f}) we calculate
\bee
\begin{split}
\hksqrt{p_f(d)+d}-\sqrt{p_f(d)}  \sim \hksqrt{\sqrt{d} e^{\sqrt{d}}+d}-\hksqrt{\sqrt{d} e^{\sqrt{d}}}=\\
\hksqrt{\sqrt{d} e^{\sqrt{d}}}\Big(\hksqrt{1+\frac{d}{\sqrt{d} e^{\sqrt{d}}}}-1\Big)=
\frac{1}{2}d^{\frac{\small 3}{\tiny 4}}e^{-\frac{1}{2}\sqrt{d}}+\ldots .
\end{split}
\label{wykladki}
\eee
Substituting here for $d$ the maximal gap $g(x)$ given by \eqref{d_max}
we obtain the approximate formula for $R(x)$:
\bee
R(x) \sim \frac{1}{2}g(x)^{3/4}e^{-\frac{1}{2}\sqrt{g(x)}}
\label{main-formula}
\eee
The comparison with real data is given in Figure \ref{fig-Andrica}.

The maximum of the function $\frac{1}{2}x^{\frac{\small 3}{\tiny 4}}e^{-\frac{1}{2}\sqrt{x}}$
is reached at $x=9$ and has the value $0.57971\ldots$. The maximal value of $A_n$ is
$0.6708735\ldots$ for $d=4$ and second value is $0.6392819\ldots$ for $d=14$.
Let us remark that $d=9$ is exactly in the middle between $4$ and 14.

Because in (\ref{main-formula}) $R(x)$ contains exponential of $\sqrt{g(x)}$, it is
very  sensitive to the form of $g(x)$.  The substitution  $g(x)=\log^2(x)$
leads to the form:
\bee
R(x)=\frac{\log^{3/2}(x)}{2\sqrt{x}}.
\label{R-Cramer}
\eee
This form of $R(x)$ is plotted in Fig.\ref{fig-Andrica} in green. If we will use the
guess   $p_f(d) \sim e^{\sqrt{d}}$ (\ref{Shanks1}) made by D. Shanks then
we will get the expression:
\bee
\sqrt{p_f(d)+d}-\sqrt{p_f(d)}=\frac{1}{2} d e^{-\frac{1}{2}\sqrt{d}}
\label{Shanks-2}
\eee
instead of (\ref{wykladki}). Substitution here for $d$ the form (\ref{d_max2})
leads to the curve plotted in Fig.\ref{fig-Andrica} in blue.

Finally, let us remark  that  from the above analysis it follows  that
\bee
\lim_{n \to \infty} (\sqrt{p_{n+1}} - \sqrt{p_n}) =0
\eee
The above  limit  was mentioned on p. 61 in \cite{Golomb1976} as a difficult problem
(yet unsolved).\\

\section{Conclusions}
\label{conclusion}

We have formulated  a few  conjectures on the gaps between consecutive primes,
in particular we have expressed maximal gap $G(x)$ directly by $\pi(x)$.
The guessed formulas are well confirmed by existing computer data.
The proofs  of them seem to be far away and in conclusion  we  quote here the
following remarks of R. Penrose from \cite{Penrose}, p.422:
\begin{quotation}{\it
Rigorous argument is usually the {\it last} step! Before that, one has to make
many guesses, and for these, aesthetic convictions are enormously important ---
always constrained by logical arguments and known facts.}
\end{quotation}

{\bf Acknowledgment:} I thank Professor  A.  Schinzel for correspondence and
Professor J.  B{\"u}the,   Professor   A. Granville, Professor  A.M. Odlyzko
for e-mail exchange and remarks and  Professor  W.Narkiewicz
and  my friend Dr J. Cis{\l}o for discussions and  comments.

\end{document}